%% file: do.tex
\documentclass[dvipsnames,onefignum,onetabnum,final]{siamart220329}
\title{Computing accurate eigenvalues using
  a mixed-precision Jacobi algorithm}
\author{Nicholas J.~Higham%
  \thanks{The author is deceased. Former address: Department of Mathematics,
    University of Manchester,
    Manchester M13 9PL,
    United~Kingdom.}
  \and Fran\c{c}oise~Tisseur\thanks{Department of Mathematics,
    University of Manchester,
    Manchester M13 9PL,
    United~Kingdom.}
  \and Marcus~Webb\footnotemark[2]
  \and Zhengbo~Zhou\footnotemark[2] \footnote{C\MakeLowercase{orresponding author (\email{zhengbo.zhou@student.manchester.ac.uk}).}}
  }

\headers{A Mixed-precision Jacobi algorithm}{
  Nicholas J.~Higham,
  Fran\c{c}oise~Tisseur,
  Marcus~Webb,
  Zhengbo~Zhou}

\usepackage{amsfonts}
\usepackage{amsmath}
\usepackage{amssymb} 
\usepackage{bookmark}
\usepackage{booktabs}
\usepackage{enumerate}
\usepackage{enumitem}
\usepackage{graphics,graphicx}
\usepackage[showonlyrefs]{mathtools}
\usepackage[utf8]{inputenc}
\usepackage{upref}
\usepackage{subdepth}
\usepackage[normalem]{ulem}
\usepackage{enumitem}
\usepackage{comment}

\makeatletter
\newcommand{\myitem}[2][]{\item[#2#1]\protected@edef\@currentlabel{#2}}
\makeatother

\usepackage{tcolorbox}
\tcbuselibrary{listings,breakable}
\colorlet{LightCornflowerBlue}{CornflowerBlue!50}
\newtcolorbox{mybox}{
  colback=LightCornflowerBlue,
  colframe=black,
  boxrule=0pt,
  breakable,
  arc=0pt,
  outer arc=0pt
}

\usepackage{caption}
\usepackage{subfig}
\usepackage{tikz,pgfplots,pgfplotstable,pgf}
\usepackage{subfig}
\usetikzlibrary{external}
\usepgfplotslibrary{groupplots}
\pgfplotsset{compat=1.18}

\pgfplotsset{
  compat=newest, 
  width=0.45\textwidth,
  height=0.36\textwidth,
  grid=both,
  major grid style={dashed, gray},
  minor grid style={dotted, lightgray},
  enlargelimits = false,
  axis background/.style={fill=white},
  every axis plot/.append style={
    line width=1.0pt
  },
  legend style={
    anchor=center,
    draw=black,
    legend cell align={left},
    font=\scriptsize
  }
}

\makeatletter
\g@addto@macro\bfseries{\boldmath}
\makeatother

\DeclareSymbolFont{CMlargesymbols}{OMX}{cmex}{m}{n}
\DeclareMathDelimiter{(}{\mathopen} {operators}{"28}{CMlargesymbols}{"00}
\DeclareMathDelimiter{)}{\mathclose}{operators}{"29}{CMlargesymbols}{"01}
\DeclareMathAlphabet\mathcal{OMS}{cmsy}{m}{n}
\SetMathAlphabet\mathcal{bold}{OMS}{cmsy}{b}{n}

\usepackage[lined,algo2e,linesnumbered,ruled]{algorithm2e}
\DontPrintSemicolon
\SetKwInput{kwin}{Input}
\SetKwInput{kwout}{Output}

\newsiamthm{assumption}{Assumption}
\newsiamremark{remark}{Remark}

\usepackage[hyperpageref]{backref}
\renewcommand*{\backref}[1]{}
\renewcommand*{\backrefalt}[4]{%
\ifcase #1 %
\or
   (Cited on p.~#2)%
\else
   (Cited on pp.~#2)%
\fi
}

\def\ycite[#1#2#3#4#5]#6{\cite[$\mit{#1#2#3#4}$#5]{#6}}

\mathchardef\Gamma="7100
\mathchardef\Delta="7101
\mathchardef\Theta="7102
\mathchardef\Lambda="7103
\mathchardef\Xi="7104
\mathchardef\Pi="7105
\mathchardef\Sigma="7106
\mathchardef\Upsilon="7107
\mathchardef\Phi="7108
\mathchardef\Psi="7109
\mathchardef\Omega="710A

\newcommand{\R}{\mathbb{R}}

\def\E{\mathcal{E}}

\def\eu{{\mathrm{e}}}

\newcommand{\wt}{\widetilde}

\newcommand{\n}{^{n}}

\newcommand{\nn}{^{n\times n}}
\newcommand{\tp}{^{T}}

\newcommand{\diag}{\operatorname{diag}}

\newcommand{\abs}[1]{\lvert{#1}\rvert}

\newcommand{\norm}[1]{\|{#1}\|}
\newcommand{\fnorm}[1]{\norm{#1}_F}
\newcommand{\tnorm}[1]{\norm{#1}_2}

\newcommand{\mnorm}[1]{\norm{#1}_{\mathrm{max}}}

\newcommand{\off}{{\operatorname{off}}}
\DeclareMathOperator{\fl}{\operatorname{f\kern.2ptl}}

\newcommand{\iter}[1]{^{(#1)}}
\newcommand{\tol}{{\mathrm{tol}}}
\newcommand{\eps}{{\varepsilon}}

\newcommand{\ul}{{u_{\ell}}}
\newcommand{\uh}{{u_{h}}}

\usepackage{listings}
\definecolor{gray}{rgb}{0.5,0.5,0.5}
\definecolor{mauve}{rgb}{0.58,0,0.82}
\definecolor{lightgrey}{rgb}{0.9,0.9,0.9}
\definecolor{darkgreen}{rgb}{0,0.6,0}
\lstdefinestyle{mystyle}{%
  language=matlab,
  morekeywords={anymatrix},
  showstringspaces=false,
  columns=flexible,
  keepspaces = true,
  basicstyle={\small\ttfamily},
  numbers=none,
  numberstyle=\tiny\color{gray},
  keywordstyle=\color{blue},
  commentstyle=\color{darkgreen},
  stringstyle=\color{mauve},
  breakatwhitespace=true,
  upquote=true,
  xleftmargin=5.0ex
}
\def\inline{\lstinline[basicstyle=\upshape\ttfamily]}

\def\eval{eigenvalue}

\def\evals{eigenvalues}

\def\evec{eigenvector}

\def\spd{symmetric positive definite}

\def\alg{algorithm}

\def\2nsit{2--step NS iteration}

\def\Ql{{Q_{\ell}}}

\def\At{{\wt{A}}}

\def\Qt{{\wt{Q}}}

\def\Dt{{\wt{D}}}

\newcommand{\Atcomp}{{\wt{A}_{\mathrm{comp}}}}
\newcommand{\Athcomp}{{\wt{A}_{\mathrm{h,comp}}}}

\def\ferrk{{\eps_{fwd}\iter{k}}}

\def\Qpolar{Q_{\rm polar}}
\def\QQR{Q_{\rm QR}}
\def\rfe{relative forward error}

\def\at{\wt{a}}

\def\gammah{\gamma_{h}}
\def\scond{\kappa_{2}^{S}}
\def\k{\kappa_{2}}

\crefname{algocf}{Algorithm}{Algorithms}

\begin{document}
\maketitle

\begin{abstract}
We provide a rounding error analysis of
a mixed-precision preconditioned Jacobi algorithm,
which uses low precision to compute the preconditioner, applies it at high precision
(amounting to two matrix-matrix multiplications)
and solves the eigenproblem using the Jacobi algorithm at working precision.
Our analysis yields meaningfully smaller relative forward error bounds for the computed eigenvalues compared with those of the Jacobi algorithm.
We further prove that, after preconditioning, if the off-diagonal entries of the preconditioned matrix are sufficiently small relative to its smallest diagonal entry, the relative forward error bound is independent of the condition number of the original matrix.
We present two constructions for the preconditioner
that exploit low precision, along with their error analyses. 
Our numerical experiments confirm
our theoretical results and compare the relative forward error of the proposed algorithm with the Jacobi algorithm,
a preconditioned Jacobi algorithm, and MATLAB's \texttt{eig} function. Timings using Julia suggest that the dominant cost of obtaining this level of accuracy comes from the high precision matrix-matrix multiplies; if support in software or hardware for this were improved, then this would become a negligible cost. 
\end{abstract}



\section{Introduction}
\label{sec.introduction}
The Jacobi (eigenvalue) algorithm~\cite{jaco46} is an iterative method to compute the spectral decomposition of a symmetric matrix.
Demmel and Veseli\'{c}~\cite{deve92} showed that 
for a \spd\ matrix $A\in\R\nn$, 
the eigenvalues computed by the Jacobi algorithm
with the stopping criterion
\begin{equation}\label{deve-stopping-criterion}
  \text{stop when }
  \abs{a_{ij}} \le \tol \sqrt{a_{ii}a_{jj}}
  \text{ for all $i,j$,}
\end{equation}
where $\tol > 0$ is usually a multiple of 
the working precision $u$, satisfy
\begin{equation}\label{eq.deve92-bound}
  \frac{\abs{\widehat{\lambda}_{k}(A) - \lambda_{k}(A)}}{
  \abs{\lambda_{k}(A)}} \le 
  p(n)u\scond(A),
\end{equation}
where $\lambda_{k}(A)$ and $\widehat{\lambda}_{k}(A)$ denote the $k$th largest \eval\ and computed eigenvalue of $A$,
$p(n)$ is a low-degree polynomial, and $\scond(A)$ is the
\emph{scaled condition number} of $A$ defined by
\begin{equation}\notag
  \scond(A)= \k(DAD), \quad
  D = \diag(a_{ii}^{-1/2}),
\end{equation}
with $\k(A)=\lambda_{1}(A)/\lambda_{n}(A)$.
The right-hand side of~\eqref{eq.deve92-bound}
provides a stronger bound than the relative forward error bounds
for methods that involve a tridiagonalization,
such as the divide and conquer \alg\ and the tridiagonal QR \alg,
whose bounds are proportional to $\k(A)$~\cite[sec.~4.7.1]{abbb99-LUG3}.
For example, Mathias~\cite[p.~979]{math95}
commented that if $A$ is graded then $\scond(A) \ll \k(A)$.

Most analyses of Jacobi's algorithm involve the so-called \emph{off-quantity}, $\off(A) = (\sum_{i\neq j} a_{ij}^{2})^{1/2} =\fnorm{A-\diag(A)}$. It is known that $\off(A)$ converges to zero at least linearly, and when $\off(A)$ is sufficiently small, the convergence rate is quadratic (see, for instance,~\cite{scho64,vank66,hari91}). 
From this follows the idea of \emph{preconditioning} for the Jacobi algorithm: 
construct an orthogonal matrix $\Qt$ such that $\off(\Qt\tp\!A\Qt) \ll \off(A)$
and apply
the Jacobi algorithm to $\At := \Qt\tp\!A\Qt$.
Drma\v{c} and Veseli\'{c}~\cite{drve00} developed a procedure
that uses an approximate eigenvector matrix $\Qt$ of $A$ as a preconditioner.
They also discussed a one-sided Jacobi algorithm preconditioned by one step of the Rutishauser LR diagonalization~\cite{drve08i,drve08ii}.
These preconditioning ideas for the Jacobi algorithm focus on speeding up the algorithm while maintaining the accuracy of the computed eigenvalues.
In contrast, we focus on how to use a preconditioner to \emph{improve the accuracy of the eigenvalues} computed by the Jacobi algorithm. 
The application that motivates us, but is not developed in this paper, is to improve the accuracy of spectral decomposition-based algorithms for matrix functions based on ~\cite{high08-FM}.

In recent years, researchers have sought to use the accuracy of high-precision floating point arithmetic to overcome the instability of numerical algorithms~\cite{hima22}.
Veseli\'{c} and Hari~\cite{veha89} proposed a two-stage algorithm for the one-sided Jacobi algorithm in which one performs Cholesky decomposition at high precision as a preprocessing step, then applies the one-sided Jacobi algorithm to the computed Cholesky factor. For ill-conditioned matrices, their numerical experiments yielded computed eigenvalues that were more accurate than those of the standard one-sided Jacobi algorithm.
Malyshev~\cite{maly92} and Drygalla~\cite{dryg08} suggested that applying the preconditioner to $A$ at higher precision can improve the accuracy of the computed eigenvalues. However, Malyshev only discussed the backward error, and Drygalla observed the high accuracy property without proof.

Our first contribution is to provide a theoretical justification for Drygalla's observation.
Under reasonable assumptions discussed in section \ref{sec.main-result-proof}, we prove that applying the preconditioner to $A$ at higher precision followed by the Jacobi algorithm at working precision yields computed eigenvalues that have relative forward errors bounded by $p(n)u\scond(\At)$ (see Theorem~\ref{thm.main}).
Compare this with the error bound in~\eqref{eq.deve92-bound}. 
In addition, we establish bounds on the magnitude of $\scond(\At)$ in terms of $\off(\At)$. For instance, we prove that if $\off(\At)/\min_i(\at_{ii}) < 1/2$,
then $\scond(\At) < 3$.

Our second contribution relates to constructing preconditioners for the Jacobi algorithm using a lower precision. We discuss two approaches. The first approach is to construct the preconditioner $\Qt$ by a low-precision spectral decomposition and then orthogonalize $\Qt$ to working precision (see \cite{zhba22} and \cite{zhou22-MPJA}). 
In the second approach, we perform the tridiagonalization step at a lower precision and then a tridiagonal eigensolver at the working precision. The key insight is that the low-precision orthogonal transformations arising from the tridiagonalization step can be made orthogonal to working precision almost for free, by exploiting the nature of Householder reflectors.
We prove that both approaches
construct a preconditioner $\Qt$ satisfying $\tnorm{\Qt\tp \Qt - I} \lesssim u$, and produce a preconditioned matrix $\At$ satisfying
$\off(\At)/\min_i(\at_{ii}) \lesssim \ul\kappa_{2}(A)$,
both of which are relevant to the forward error analysis in the main theorem of the paper.

The results presented in this paper have the potential to be developed into a practical algorithm for computing eigenvalues of symmetric positive definite matrices to high relative accuracy, in which the lower and higher precisions are selected based on an estimate for $\kappa_2(A)$. 

The rest of the article is organized as follows.
In section~\ref{sec.main-result-proof},
we present the
mixed-precision preconditioned Jacobi algorithm
and then state and prove Theorem~\ref{thm.main}, our main result.
In the rest of section~\ref{sec.main-result-proof},
we discuss the selection of $\uh$, and compare $\scond(\At)$ and $\scond(A)$.
In section~\ref{sec.altern-constr-prec} we provide and analyze methods for
constructing an approximate eigenvector matrix
as a preconditioner.
Numerical results are presented in section~\ref{sec.numerical-examples}
to support our theoretical analysis.
Finally, conclusions and further remarks are given in
section~\ref{sec.conclusion}.


\section{Mixed-precision preconditioned Jacobi algorithm}
\label{sec.main-result-proof}

We are interested in the computation of a spectral decomposition of 
a symmetric matrix $A\in\R\nn$
and its \rfe, defined as
\begin{equation}
\label{def-ferrk-specific}
  \ferrk =
  \frac{\abs{\widehat{\lambda}_{k}(\Atcomp) -
    \lambda_{k}(A)}}{\lambda_{k}(A)},
\end{equation}
where $\Atcomp$ is the computed preconditioned matrix $\At=\Qt\tp A\Qt$. 
Throughout, we treat $u$ and $\uh$ with $u>\uh$
as working and high precision, respectively.
They are chosen from floating-point formats defined by IEEE standard~\cite[Tab.~3.5]{ieee19} and are assumed to satisfy $n\uh < nu < 1/2$ and $u < 0.1$.


Demoting a floating-point number from one precision to another that is lower than the current one can lead to overflow, underflow, or subnormal numbers. 
For Algorithms~\ref{alg.orth-preconditioner} and~\ref{alg.modified-tridiag}, which are described in section~\ref{sec.altern-constr-prec} for constructing the preconditioner $\Qt$, underflow is not detrimental since the relevant analysis is purely normwise. 
On the other hand, for Algorithm~\ref{alg.uh-prec-jacobi}, the occurrence of underflow compromises the componentwise analysis. So, throughout this work, we assume that underflow and overflow do not occur. 


The following is the main algorithm of the paper, the mixed-precision preconditioned Jacobi algorithm, which uses a high precision $\uh$ to apply the preconditioner.
Throughout the paper, $p_k \coloneq p_k(n)$ is a low-degree polynomial.

\begin{algorithm2e}[h]
\caption{Mixed-precision preconditioned Jacobi algorithm.
  \label{alg.uh-prec-jacobi}}
\kwin{A symmetric matrix $A\in\R\nn$,
  two precisions $u$ and $\uh$ with $0 < \uh < u$,
  and a preconditioner $\Qt$ satisfying
  $\tnorm{\Qt\tp\Qt-I} \le p_1u < 1$.}
\kwout{A computed spectral decomposition $Q\Lambda Q\tp$
  of $A$, where $\Lambda\in\R\nn$ is diagonal and
  $\tnorm{Q\tp Q - I} =O(u)$.}
\BlankLine
Compute the preconditioned matrix $\At$
by computing the product $\Qt\tp A \Qt$ entirely at precision $\uh$,
which gives $\Athcomp$, and then demoting it to precision $u$ to obtain $\Atcomp$. \;\label{li.uh-pjb.use-of-uh}
Compute a spectral decomposition
$Q_{J} \Lambda Q_{J}\tp$ of $\Atcomp$
using the Jacobi algorithm with
stopping criterion~\eqref{deve-stopping-criterion}
at precision $u$.\;\label{li.uh-prec-jacobi:Jacobi}
Construct the eigenvector matrix $Q = \Qt Q_{J}$ at precision $u$.\;\label{li.Form-the-eigenvector-matrix}
\end{algorithm2e}

In line~\ref{li.uh-pjb.use-of-uh}, to ensure that the computed $\Atcomp$ is symmetric, one can compute its symmetric part via $(\Atcomp+\Atcomp\tp)/2$ or make use of BLAS 3 routines as described in \cite[Rmk.~2.2]{drve00}.

Efficient constructions of such $\Qt$  are discussed in
section~\ref{sec.altern-constr-prec}. 
Although Algorithm~\ref{alg.uh-prec-jacobi} is able to compute the full spectral decomposition of a symmetric matrix, the main theorem (Theorem \ref{thm.main}) applies only to symmetric positive definite $A$, so we assume positive definiteness for the remainder of the section.

\subsection{The main theorem}
\label{sec.setting-scene}

Following \cite{high02-ASNA2}, define 
\begin{equation}\label{gamma-gammah}
  \gamma_{h} \coloneq \frac{nu_{h}}{1-nu_{h}} <
  \gamma \coloneq \frac{nu}{1-nu} < 1.
\end{equation}
Our main theorem (Theorem \ref{thm.main}) uses the following assumptions.
\begin{assumption}
\label{ass.main}
Let $u,\uh,\gamma$ and $\gammah$ be as in \eqref{gamma-gammah}, and let 
$A, \At \in \R\nn$ and $p_1$ be as in
Algorithm~\ref{alg.uh-prec-jacobi}. We assume
\begin{enumerate}[label=\textup{(A\arabic*)}]
\item \label{it.ass.nukappa}
$10n^{3/2}u (1-p_1u)^{-1} \kappa_2(A) < 1$,
\item \label{it.ass-size-gammah}
$\gammah <  \frac{u(1-p_1u)}{16n^{1/2}\kappa_2(A)}$,
\item \label{it.c(n,u)kappa<1/2}
$14nu\scond(\At) < 1$.
\end{enumerate}
\end{assumption}

Notice that the assumption that $p_1u < 1$ in
Algorithm~\ref{alg.uh-prec-jacobi}
is necessary for Assumptions \ref{it.ass.nukappa} and \ref{it.ass-size-gammah} to be satisfiable.  We use Assumption \ref{it.ass.nukappa} to guarantee that $\Atcomp$ is positive definite, and to control a term of the form $O(u^2 \kappa_2(A))$ in the error bounds.  We use Assumption \ref{it.c(n,u)kappa<1/2} to bound $\scond(\Atcomp)$ in terms of $\scond(\At)$. These assumptions are not strictly necessary to obtain an informative bound on $\ferrk$; indeed, we derive several intermediate results based on fewer assumptions during the proof of Theorem \ref{thm.main}.

\begin{theorem}
\label{thm.main}
Let $A\in\R\nn$ be symmetric positive definite, and let
$\At =\Qt\tp\!A\Qt$ be the preconditioned matrix as in Algorithm~\ref{alg.uh-prec-jacobi}.
If Assumption~\ref{ass.main} holds, then
the forward error $\ferrk$ in~\eqref{def-ferrk-specific} satisfies
\begin{equation}\label{eq:ferrk-main}
  \ferrk \le p_3 \scond(\At) u + p_1u,
\end{equation}
provided that both terms on the right-hand side are less than $1$.
Here $p_3$ is a polynomial in $n$ that is at most quadratic.
\end{theorem}

It is possible to absorb the second term on the right-hand side of~\eqref{eq:ferrk-main} into the first. However, we retain the splitting because the first term describes the error introduced by applying the preconditioner and the eigensolver,
while the second term captures the error that arises
during the construction of the preconditioner.



We divide the proof of Theorem~\ref{thm.main} into three parts
corresponding to the decomposition,
\begin{align}
  \ferrk
  & \le
    \frac{\abs{\widehat{\lambda}_{k}(\Atcomp) -
    \lambda_{k}(\Atcomp)}}{\lambda_{k}(A)} +
    \frac{\abs{\lambda_{k}(\Atcomp) -
    \lambda_{k}(\wt A)}}{\lambda_{k}(A)} +
    \frac{\abs{\lambda_{k}(\wt A) -
    \lambda_{k}(A)}}{\lambda_{k}(A)}\notag \\
  & \eqcolon E_{1} + E_{2} + E_{3}. \label{eq.proof-decomp}
\end{align}
Bounding $E_3$ is straightforward.  
Ostrowski's theorem \cite{ostrowski1959quantitative} implies that
$\lambda_{k}(\At) = \theta_{k} \lambda_{k}(A)$
with $\lambda_n(\Qt\tp\Qt) \le\theta_{k} \le \lambda_{1}(\Qt\tp\Qt)$.
Since $\Qt$ is assumed to satisfy $\fnorm{\Qt\tp\Qt-I} \le p_1u$,
namely $\Qt\tp\Qt = I + E$ where $\tnorm{E}\le p_1u$,
by~\cite[Cor.~8.1.6]{gova13-MC4} we have 
\begin{equation}\label{lower-bound-for-lambda}
    |\lambda_n(\Qt\tp\Qt)-1|
    = |\lambda_n(I+E)-\lambda_n(I)| \le p_1u.
\end{equation}
Thus, $1-p_1u\le\lambda_n(\Qt\tp\Qt)$.
Similarly, we get $\lambda_1(\Qt\tp\Qt) \le 1+p_1u$.
Then, $-p_1u \le \theta_k-1 \le p_1u$, 
namely $|\theta_k - 1| \le p_1u$.
Therefore, we can bound $E_3$ as follows: 
\begin{equation}
    \label{err-E3}
    E_3 = \frac{|\theta_k\lambda_k(A)-\lambda_k(A)|}{\lambda_k(A)}
    = |\theta_k-1| \le p_1 u.
\end{equation}

Before bounding $E_1$ and $E_2$, we derive a useful bound for $|\Delta \wt A|$, 
the elementwise absolute value of the error $\Delta \wt A=\Atcomp - \wt A$. 
\begin{lemma}\label{lem.At-Atcomp}
Let $\At, \Atcomp \in \R\nn$ be  as in Algorithm~\ref{alg.uh-prec-jacobi}.
Then
\begin{equation}\label{err-product}
  \Atcomp = \wt A + \Delta \wt A, \qquad
  \abs{\Delta \wt A} \le
  u\abs{\At} + 4\gammah
  \abs{\wt Q\tp}\, \abs{A}\, \abs{\wt Q}.
\end{equation}
\end{lemma}
\begin{proof}
Storage in IEEE standard formats ensures that
data promoted from working to high precision stays the same~\cite[p.~39]{over01-ncieeefpa}.
Consequently, the perturbation $\Delta\At$ comes solely from the following two operations:
\begin{enumerate}[label=\upshape(\arabic*),labelwidth=!]
\item
performing matrix products with high precision $\uh$, and
\item
demoting a matrix stored at precision $\uh$ back to the working precision $u$.
\end{enumerate}
We treat the computation of $\wt A$ in floating point arithmetic
at precision $\uh$ as
\begin{equation}
  \Athcomp = (\wt Q\tp A + F_{1})\wt Q + F_{2}
   = \Qt\tp A\Qt + (F_{1} \Qt + F_{2}),
\end{equation}
where $F_{1}$ and $F_{2}$ are perturbations
from computing the matrix products
$\Qt\tp A$ and $(\wt Q\tp A + F_{1})\wt Q$.
Standard rounding error analysis for
matrix multiplication~\cite[sec.~3.5]{high02-ASNA2}
gives the following componentwise bounds.
\begin{align*}
  \abs{F_{1}}
  & \le \gamma_{h} \abs{\wt Q\tp}\, \abs{A}, \\
  \abs{F_{2}}
  & \le \gamma_{h} \abs{\wt Q \tp A + F_{1}}\, \abs{\wt Q}
  \le \gamma_{h} \abs{\wt Q\tp A}\,\abs{\wt Q} +
  \gamma_{h} \abs{F_{1}}\,\abs{\wt Q}. 
\end{align*}
Letting $F_{3} = F_{1}\Qt + F_{2}$,
we obtain $\Athcomp = \wt A + F_{3}$, with
\begin{equation}\label{err-mmm}
  \begin{aligned}
    \abs{F_{3}}
    & \le \abs{F_{1}}\,\abs{\Qt}+\abs{F_{2}} \\
    & \le \gammah \abs{\Qt\tp}\,\abs{A}\,\abs{\Qt} +
      \gammah\abs{\Qt\tp}\,\abs{A}\,\abs{\Qt} +
      \gammah^{2}\abs{\Qt\tp}\,\abs{A}\,\abs{\Qt} \\
    & \le 3 \gammah \abs{\Qt\tp}\,\abs{A}\,\abs{\Qt}.
  \end{aligned}
\end{equation}
The last inequality comes from $\gammah^{2} < \gammah$. 


Demoting $\Athcomp$, stored in high precision, to working precision is the same as rounding $\Athcomp$ to working precision.
Using~\cite[Thm.~2.2]{high02-ASNA2}, the rounded matrix $\Atcomp$ satisfies
\begin{equation}\notag
  \Atcomp = \Athcomp + F_{4}, \quad
  \abs{F_{4}} \le u \abs{\Athcomp}. 
\end{equation}
Substituting~\eqref{err-mmm} into the inequality above yields
\begin{equation}\label{err-rounding}
  \abs{F_{4}} \le u \abs{\At} + u \abs{F_{3}} \le
  u \abs{\At} + 3u\gammah
  \abs{\wt Q\tp}\, \abs{A}\, \abs{\wt Q}.
\end{equation}
Combining~\eqref{err-mmm} and~\eqref{err-rounding}, we obtain
\begin{align}
  \abs{\Atcomp - \wt A}
  & \le \abs{\Atcomp- \Athcomp} +
    \abs{\Athcomp - \wt A} \notag\\
  & \le \abs{F_{4}} + \abs{F_{3}} \notag \\
  & \le 3 \gammah \abs{\Qt\tp}\abs{A} \abs{\Qt}
     + u \abs{\At} + 3u\gammah
      \abs{\wt Q\tp}\, \abs{A}\, \abs{\wt Q}. \label{eq.opjsd}
\end{align}
The upper bound for $\abs{\Delta \wt A}$ in \eqref{err-product}
follows since by assumption, $u<0.1$ so that $3u\gammah < \gammah$.
\end{proof}
\begin{remark}
\label{rmk.compwise-error}
The error bound~\eqref{err-product} is subtle but is key to understanding why only the matrix-matrix multiplications are required at high precision. If all computations are performed at working precision, then we obtain the bound,
\begin{equation}\label{normal-bound}
  \abs{\Delta \wt A} \le
  \gamma \abs{\wt Q\tp} \abs{A} \abs{\Qt}.
\end{equation}
The term $\abs{\wt Q\tp} \abs{A} \abs{\Qt}$ is potentially large relative to $\sqrt{\at_{ii}\at_{jj}}$, which impedes the use of perturbation bounds, such as Theorem \ref{thm.perturb-evals}. In the bound~\eqref{err-product}, we can control this term using a sufficiently high precision $u_h$, but in the bound~\eqref{normal-bound}, we cannot.
\end{remark}
\subsection{Bound on the forward error $E_1$}
\label{sec.forw-error-comp-1}
\begin{lemma}\label{lem.Atcomp-spd}
Let $\Atcomp \in \R\nn$ be as in Algorithm~\ref{alg.uh-prec-jacobi}. If Assumption~\ref{it.ass.nukappa} holds, then $\Atcomp$ is positive definite.
\end{lemma}
\begin{proof}
From the componentwise bound~\eqref{err-product} we have%
\begin{equation}\notag
  \fnorm{\Delta \wt A}
  \le u \fnorm{\At} + 4\gammah\fnorm{\Qt}^{2}\fnorm{A}
  \le (u + 4\gammah)\fnorm{\Qt}^{2}\fnorm{A}
  \le n(u+4\gammah)\tnorm{\Qt}^2\fnorm{A}.
\end{equation}
Using $\tnorm{\Qt}^2=\tnorm{\Qt\tp\Qt}\le 1+p_1u$ as well as $\gammah < u$ and $p_1u < 1$, we have
\begin{equation}\label{err-product-normwise}
  \fnorm{\Delta \wt A}
  \le n(u + 4\gammah)(1+p_1 u)\fnorm{A}
    \le n(2u + 8\gammah)\fnorm{A}
  \le 10nu\fnorm{A}.
\end{equation}

Applying the eigenvalue perturbation bound~\cite[Cor.~8.1.6]{gova13-MC4} to $\Atcomp$ together with~\eqref{err-product-normwise} yields a lower bound on $\lambda_{n}(\Atcomp)$,
\begin{equation}\notag
  \lambda_{n}(\Atcomp)
  \ge \lambda_{n}(\At) - \tnorm{\Delta\At}
  \ge \lambda_{n}(\At) - \fnorm{\Delta \At}
  \ge \lambda_{n}(\At) - 10nu\fnorm{A}.
\end{equation}
Using the fact that $\fnorm{A} \le n^{1/2}\tnorm{A}$ 
for $A\in\R^{n\times n}$,
Ostrowski's theorem~\cite{ostrowski1959quantitative},
$1-p_1u\le\lambda_n(\Qt\tp\Qt)$ (by ~\eqref{lower-bound-for-lambda}), 
and Assumption~\ref{it.ass.nukappa}, we have
\begin{equation}
    \notag
    \frac{10nu\fnorm{A}}{\lambda_n(\At)} \le 
    \frac{10n^{3/2}u \lambda_1(A)}{(1-p_1 u)\lambda_n(A)} < 1.
\end{equation}  
Therefore $10nu\fnorm{A} < \lambda_n(\At)$, and $\Atcomp$ is positive definite.     
\end{proof}

As a consequence of Lemma~\ref{lem.Atcomp-spd}, we can use 
the relative perturbation result for
applying the Jacobi algorithm on $\Atcomp$%
~\cite[Thms.~2.3 \&~3.1]{deve92}, that is,
\begin{equation}\label{fwdbd-Atcomp}
  \frac{\abs{\widehat{\lambda}_{k}(\Atcomp) -
  \lambda_{k}(\Atcomp)}}{\lambda_{k}(\Atcomp)}
  \le p_{2} u \scond(\Atcomp), 
\end{equation}
provided that $ p_{2} u \scond(\Atcomp)<1$.
In the worst-case scenario, $p_{2}$ is quadratic in $n$,
but in practice, $p_{2} \approx cn$, 
as observed in~\cite{deve92}.
Substituting~\eqref{fwdbd-Atcomp} into
$E_{1}$ in~\eqref{eq.proof-decomp} gives%
\begin{equation}\label{eq.E1.inter1}
  E_{1} =
  \frac{\abs{\widehat{\lambda}_{k}(\Atcomp) -
    \lambda_{k}(\Atcomp)}}{\lambda_{k}(A)}
  \le p_{2} u \scond(\Atcomp)
  \frac{\lambda_{k}(\Atcomp)}{\lambda_{k}(A)}. 
\end{equation}
Writing $\Atcomp$ as in~\eqref{err-product} and
using the Courant--Fischer theorem~\cite[Thm.~8.1.5]{gova13-MC4} yields
\begin{equation}\label{bound-intermediate}
  \frac{\lambda_{k}(\Atcomp)}{\lambda_{k}(A)}
  \le
    \frac{\lambda_{k}(\wt A) + \tnorm{\Delta \wt A}}{\lambda_{k}(A)}
  \le
  \frac{\lambda_{k}(\wt A) + \fnorm{\Delta \wt A}}{\lambda_{k}(A)}
  = \frac{\lambda_{k}(\wt A)}{\lambda_{k}(A)} +
  \frac{\fnorm{\Delta \wt A}}{\lambda_{k}(A)}.
\end{equation}
Using Ostrowski's theorem \cite{ostrowski1959quantitative},
$\tnorm{\Qt\tp\Qt}\le 1+p_1u$,
and the normwise bound~\eqref{err-product-normwise},
we have%
\begin{equation}\notag
  \frac{\lambda_{k}(\Atcomp)}{\lambda_{k}(A)}
  \le \tnorm{\Qt\tp\Qt} + \frac{10nu\fnorm{A}}{\lambda_{k}(A)}
  \le 1 + p_{1}u + 10n^{3/2}u\kappa_{2}(A). 
\end{equation}
Substituting this into \eqref{eq.E1.inter1} gives 
\begin{align*}
  E_{1}
  & \le p_{2}u\scond(\Atcomp)
  (1 + p_{1}u + 10n^{3/2}u\kappa_{2}(A))
  \\
  & \le \big(
  p_{2}u+p_{1}p_{2}u^{2}+10n^{3/2}p_{2}u^{2}\kappa_{2}(A)
  \big)\scond(\Atcomp).
\end{align*}
Now, $p_{1}u<1$ implies $p_{1}^{}p_{2}^{}u_{}^{2}<p_{2}^{}u$
and Assumption~\ref{it.ass.nukappa} implies $10n^{3/2}u\kappa_{2}(A)<1$. Therefore
\begin{align*}
  E_{1}
  & \le \big(2p_{2}u + 10n^{3/2}p_{2}u^{2}\kappa_{2}(A)\big)
    \scond(\Atcomp).\\
  & \le (2p_{2}u + p_{2}u)\scond(\Atcomp)\\
  & \le 3p_{2}u\scond(\Atcomp). 
\end{align*}
Let us redefine $p_{2}$ as $3p_{2}$.
We have
\begin{equation}\label{err-E1}
  E_{1} \le p_{2}u\scond(\Atcomp),
\end{equation}
where $p_{2}$ is a polynomial in $n$ that is at most quadratic. 

\subsection{Bound on the forward error $E_2$}
\label{sec.forw-error-multiplicative}
The next result is a direct consequence of~\cite[Thms.~2.1 \&~2.17]{vesl93}. See also \cite{ipsen1998relative} for related bounds.
\begin{theorem}[{\cite{vesl93}}]
\label{thm.perturb-evals}
Let $A\in\R\nn$ be symmetric positive definite and $\Delta A$ be a symmetric perturbation that satisfies
\begin{equation}\label{eq.delta-A-perturb}
  \abs{\Delta a_{ij}} \le \eps \sqrt{a_{ii}a_{jj}}
\end{equation}
for some $\eps > 0$.
If $\eps n \scond(A) < 1$,
then
\begin{equation}\label{eq.rel-eig-bound-1}
  \frac{ \abs{\lambda_{k}(A+ \Delta A)-\lambda_{k}(A)}}{\lambda_{k}(A)}
  \le \eps n \scond(A)
\end{equation}
for all $k$.
\end{theorem}
The bound~\eqref{eq.rel-eig-bound-1} is similar to~\eqref{eq.deve92-bound},
since the Jacobi algorithm with stopping criterion~\eqref{deve-stopping-criterion}
computes the true eigenvalues of $A+\Delta A$,
where $\Delta A$ satisfies~\eqref{eq.delta-A-perturb}%
~\cite[Thms. 2.3 \& 3.1]{deve92}.
Rewriting the componentwise bound~\eqref{err-product} as
\begin{equation}\label{eq.error-atij-deltaatij}
  \abs{\Delta \at_{ij}} \le
  u \abs{\at_{ij}} +
  4\gammah \big( \abs{\wt Q\tp} \abs{A} \abs{\Qt} \big)_{ij},
\end{equation}
and using the fact that $\At$ is positive definite so that $\at_{ij} \leq \sqrt{\at_{ii}\at_{jj}}$,
we have 
\begin{equation}\notag
  \abs{\Delta \at_{ij}} \le
  u \sqrt{\at_{ii}\at_{jj}} +
  4\gammah \big( \abs{\wt Q\tp} \abs{A} \abs{\Qt} \big)_{ij}.
\end{equation}
Rearranging the inequality yields
\begin{equation}\label{eq.high-relative-accuracy}
  \frac{\abs{\Delta \at_{ij}}}{\sqrt{\at_{ii}\at_{jj}} }
  \le u + 4 \gammah \alpha_{ij},\qquad
  \alpha_{ij} \coloneq
  \frac{\big( \abs{\wt Q\tp} \abs{A} \abs{\Qt}
  \big)_{ij}}{\sqrt{\at_{ii}\at_{jj}} }.
\end{equation}
Let 
\begin{equation}\label{eq:definition-alpha}
\alpha = \max_{i,j}\alpha_{ij}.
\end{equation}
Applying Theorem~\ref{thm.perturb-evals} with $\eps = u+4\gammah\alpha$ leads to 
\begin{equation}\label{fwdbd-Atcomp-2}
  \frac{\abs{\lambda_{k}(\Atcomp)-\lambda_{k}(\At)}}{\lambda_{k}(\At)}
  \le n(u+4\gammah\alpha)\scond(\At),
\end{equation}
provided that $n(u+4\gammah\alpha)\scond(\At)<1$.
Substituting~\eqref{fwdbd-Atcomp-2} into $E_{2}$
in~\eqref{eq.proof-decomp} gives
\begin{equation}\notag
  E_{2} =
  \frac{\abs{\lambda_{k}(\Atcomp)-\lambda_{k}(\At)}}{\lambda_{k}(A)}
  \le n(u+4\gammah\alpha)\scond(\At)
  \frac{\lambda_{k}(\At)}{\lambda_{k}(A)}. 
\end{equation}
By Ostrowski's Theorem \cite{ostrowski1959quantitative}, $\lambda_{k}(\At)/\lambda_{k}(A) \le 1+p_{1}u < 2$.
Therefore, the inequality becomes%
\begin{equation}\label{err-E2}
  E_{2}
  \le \big(2nu + 8n\gammah\alpha\big)
  \scond(\At). 
\end{equation}

\subsection{An intermediate bound}
\label{sec.disc-forw-error}
Summing the bounds on $E_1$, $E_2$, and $E_3$ in ~\eqref{err-E1}, \eqref{err-E2}, and~\eqref{err-E3} leads to the following intermediate theorem.
\begin{theorem}
\label{thm.original-main}
Let $\wt A, \Atcomp,\Qt \in\R\nn$ be
as in Algorithm~\ref{alg.uh-prec-jacobi}.
If Assumption~\ref{it.ass.nukappa} and $\gammah<u$ hold,
then the forward error $\ferrk$ in~\eqref{def-ferrk-specific}
satisfies %
\begin{equation}\label{fwd-err-grand}
  \ferrk \le p_{2}u\scond(\Atcomp)
  + \big(2nu + 8n\gammah \alpha \big)
  \scond(\At) + p_{1}u,
\end{equation}
provided that each of the three terms on the right-hand side is less than $1$.
Here $p_{2}$ is a polynomial in $n$ that is at most quadratic, and
$\alpha$ is given in~\eqref{eq:definition-alpha}.
\end{theorem}


\subsubsection{On the size of $\alpha$ in \eqref{eq:definition-alpha}}
\label{sec.disc-size-alpha}
We use the next result directly to bound $\alpha$.

\begin{theorem}
[{\cite[Thm.~1.1]{chu95}}]
\label{thm.majorization}
Let $\lambda_n(A)$ and $d_n(A)$ denote the smallest
eigenvalue and diagonal entry of a symmetric $A\in\R\nn$.
Then $\lambda_n(A) \le d_n(A)$.
\end{theorem}
We can now bound $\alpha$ by 
\begin{equation}\notag
  \alpha = \max_{i,j}\alpha_{ij} \le
  \frac{\max_{i,j}(\abs{\wt Q\tp}\,\abs{A}\,\abs{\wt Q})_{ij}}{
  \min_{i}(\at_{ii})}
  \le
  \frac{\max_{i,j}(\abs{\wt Q\tp}\,\abs{A}\,\abs{\wt
  Q})_{ij}}{\lambda_{n}(\wt A)},
\end{equation}
where we use Theorem~\ref{thm.majorization} for the last inequality.
The numerator is 
$\mnorm{\abs{\wt Q\tp}\,\abs{A}\,\abs{\wt Q}}$,
where $\mnorm{\,\cdot\,}$ is the matrix max norm.
By definition~\cite[p.~110]{high02-ASNA2},
\begin{equation}\notag
  \mnorm{\abs{\wt Q\tp}\,\abs{A}\,\abs{\wt Q}}
  = \norm{\abs{\wt Q\tp}\,\abs{A}\,\abs{\wt Q}}_{1,\infty}
  \le
  \norm{\abs{\wt Q \tp}}_{2,\infty}
  \tnorm{\abs{A}}
  \norm{\abs{\wt Q}}_{1,2},
\end{equation}
where for $C \in \R\nn$,
$\norm{C}_{2,\infty} = \max_{i} \tnorm{C(i,:)\tp}$
and $\norm{C}_{1,2} = \max_{j} \tnorm{C(:,j)}$.
Therefore, we have 
$\norm{\abs{\wt Q \tp}}_{2,\infty} = \norm{\abs{\wt Q}}_{1,2}
= \norm{\wt Q}_{1,2}$.
The last equality follows from the definition of the vector $2$-norm,
which eliminates the effect of taking the absolute value.
In conclusion, we have
\begin{equation}\notag
  \alpha \le
  \frac{\norm{\wt Q}_{1,2}^{2} \tnorm{\abs{A}}}{\lambda_{n}(\wt A)}
  \le
  \frac{(1+p_{1}u)^{2} n^{1/2} \tnorm{A}}{(1-p_{1}u) \lambda_{n}(A)}
  = \frac{(1+p_{1}u)^{2} n^{1/2}\kappa_{2}(A) }{1-p_{1}u} 
  \le \frac{4n^{1/2}\kappa_{2}(A)}{1-p_{1}u}.
\end{equation}

The next corollary shows how to choose $\uh$ such that
we can drop the dependency on $\alpha$ in the bound~\eqref{fwd-err-grand}.

\begin{corollary}
\label{cor.main-thm}
Let $A, \wt A, \Atcomp,\Qt \in\R\nn$ be
as in Algorithm~\ref{alg.uh-prec-jacobi}.
If Assumptions~\ref{it.ass.nukappa} and~\ref{it.ass-size-gammah} hold,
then the forward error $\ferrk$ in~\eqref{def-ferrk-specific} satisfies
\begin{equation}\label{corollary-grand-error-bound}
  \ferrk \le p_{2}u\scond(\Atcomp) + 4nu\scond(\At) + p_{1}u,
\end{equation}
provided that each of the three terms on the right-hand side is less than $1$.
Here $p_{2}$ is a polynomial in $n$ that is at most quadratic.
\end{corollary}

\begin{proof}
From~\eqref{fwd-err-grand},
it is sufficient to prove that $4\gamma_h\alpha \le u$. 
Using the upper bound for $\gammah$ in Assumption~\ref{it.ass-size-gammah},
we have 
\begin{equation}\notag
  4\gammah\alpha
  < 4 \biggl( \frac{u(1-p_{1}u)}{16n^{1/2}\kappa_{2}(A)} \biggr)
  \biggl( \frac{4n^{1/2}\kappa_{2}(A)}{1-p_{1}u} \biggr) = u.
\end{equation}
Substituting this into~\eqref{fwd-err-grand}
yields the desired bound. 
\end{proof}

In a practical implementation with the ability to simulate arbitrary precisions, using, for example, the MATLAB Multiprecision Computing Toolbox, we can choose $u_h$ dynamically based on
Assumption~\ref{it.ass-size-gammah}.
Suppose $A\in\R^{100\times 100}$
with $\kappa(A) = 10^8$. To compute eigenvalues with $\ferrk \lesssim u$,
instead of taking $\uh$ to be quadruple precision with $34$ digits of
accuracy, a precision with $25$ digits of accuracy is sufficient in practice.

\subsubsection{$\scond(\At)$ and $\scond(\Atcomp)$}
We provide a relationship between the two condition numbers appearing in Theorem~\ref{thm.original-main} and Corollary~\ref{cor.main-thm}.

\begin{lemma}
\label{lem.bound-scAtcomp-scAt}
Let $\At, \Atcomp\in\R\nn$ be the matrices defined in
Algorithm~\ref{alg.uh-prec-jacobi},
and let $\Delta\At$ satisfy~\eqref{eq.high-relative-accuracy}.
If $\gamma_h$ satisfies Assumption~\ref{it.ass-size-gammah},
then 
\begin{equation}\label{bound-scAtcomp-scAt}
  \scond(\Atcomp) \le \frac{1 + c(n,u) }{1 - c(n,u) \scond(\At)} \scond(\At),
\end{equation}
provided $c(n,u)\scond(\At) < 1$, where
\begin{equation}\label{def.cnu}
  c(n,u) \coloneq n
  \left( 2u + \frac{2u}{\sqrt{1-2u}} + \frac{2u}{1-2u} \right).
\end{equation}
\end{lemma}
\begin{proof}
It follows from~\eqref{eq.high-relative-accuracy} that
\begin{equation}\notag
  (\Atcomp)_{ii} = (\At + \Delta \At)_{ii}=
  \at_{ii} + \Delta\at_{ii}, \qquad
  \abs{\Delta\at_{ii}} \le
  \left( u + 4\gamma_h \alpha_{ii} \right) \abs{\at_{ii}}.
\end{equation}
Let $\uh$ be chosen such that $\gamma_h$
satisfies Assumption~\ref{it.ass-size-gammah}.
Then from the proof of Corollary~\ref{cor.main-thm}, we have that 
$4\gamma_h\max_{i}\alpha_{ii} < u$ so that
$\abs{\Delta\at_{ii}} \le 2u \abs{\at_{ii}} = 2u\at_{ii}$. Hence,
\begin{equation}\label{upper-lower-bnds}
  \sqrt{\at_{ii}} \sqrt{1-2u}
  = \sqrt{\at_{ii} - 2u \at_{ii}}
  \le \sqrt{(\Atcomp)_{ii}}
  \le \sqrt{\at_{ii} + 2u \at_{ii}}
  = \sqrt{\at_{ii}} \sqrt{1+2u}.
\end{equation}
Pre- and post-multiply 
$\Atcomp = \At + \Delta \At$ with $D_1 = \diag\big( (\Atcomp)_{ii}^{-1/2} \big)$ to obtain
\begin{equation}\label{relation-D1AD1-D2AD2}
  D_1\Atcomp D_1 = D_1\At D_1 + D_1\Delta\At D_1 = D_2\At D_2 + \mathcal{F},
\end{equation}
where 
$$\mathcal{F} = (D_1 - D_2)\At D_2 + D_1\At(D_1 - D_2) + D_1\Delta\At D_1
$$ 
with 
$D_2 = \diag (\at_{ii}^{-1/2})$.
For the $i$th diagonal entry of $D_1 - D_2$, we have that
\begin{equation}\notag
  \abs{(D_1 - D_2)_{ii}} =
  \left|\frac{1}{\sqrt{(\Atcomp)_{ii}}} -
  \frac{1}{\sqrt{\at_{ii}}} \right| =
  \frac{\left|{\sqrt{\at_{ii}} - \sqrt{(\Atcomp)_{ii}}}\right|}{
  \sqrt{(\Atcomp)_{ii}}\sqrt{\at_{ii}}}.
\end{equation}
Using~\eqref{upper-lower-bnds}, we bound the
numerator above by $\sqrt{\at_{ii}}(1-\sqrt{1-2u})$.
and the denominator below by $\at_{ii} \sqrt{1-2u}$ to give
\begin{equation}\label{Bound-of-D1-D2}
  \abs{(D_1 - D_2)_{ii}} \le
  \frac{\sqrt{\at_{ii}}
  (1-\sqrt{1-2u})}{\at_{ii}\sqrt{1-2u}}
  \le \frac{2u}{\sqrt{\at_{ii}}},
\end{equation}
where the last inequality is due to $(1-\sqrt{1-2u})/\sqrt{1-2u}<2u$ for $u < 0.1$.

Using~\eqref{Bound-of-D1-D2}, we find that 
\begin{align*}
  \abs{\mathcal{F}_{ij}}
  & \le \frac{2u \abs{\at_{ij}}}{\sqrt{\at_{ii}} \sqrt{\at_{jj}}} +
    \frac{2u \abs{\at_{ij}}}{\sqrt{(\Atcomp)_{ii}} \sqrt{\at_{jj}}} +
    \frac{\abs{\Delta\at_{ij}}}{\sqrt{(\Atcomp)_{ii}}\sqrt{(\Atcomp)_{jj}}}\\
  & \le \frac{2u \abs{\at_{ij}}}{\sqrt{\at_{ii}} \sqrt{\at_{jj}}} +
    \frac{2u \abs{\at_{ij}}}{\sqrt{1-2u}\sqrt{\at_{ii}} \sqrt{\at_{jj}}} +
    \frac{\abs{\Delta\at_{ij}}}{(1-2u)\sqrt{\at_{ii}}\sqrt{\at_{jj}}}.
\end{align*}
Since $\gamma_h$ is chosen so that
$\abs{\Delta\at_{ij}}/\sqrt{\at_{ii}\at_{jj}} \le 2u$
and since $\abs{\at_{ij}} < \sqrt{\at_{ii}\at_{jj}}$
by the positive definiteness of $\At$, we conclude that
\begin{equation}\notag
  \abs{\mathcal{F}_{ij}} \le
  2u + \frac{2u}{\sqrt{1-2u}} + \frac{2u}{1-2u}.
\end{equation}
It then follows from~\cite[Lem.~6.6]{high02-ASNA2} that 
$\tnorm{\mathcal{F}} \le c(n,u)$  
with $c(n,u)$ as in \eqref{def.cnu}.
Finally, 
writing $\scond(\Atcomp)$ in terms of $\At$
using~\eqref{relation-D1AD1-D2AD2} yields
\begin{equation}\notag
  \scond(\Atcomp)
  = \frac{\lambda_1(D_2\At D_2+\mathcal{F})}{
  \lambda_n(D_2\At D_2+\mathcal{F})} 
  \le \frac{\lambda_1(D_2\At D_2) + \tnorm{\mathcal{F}}}{
    \lambda_n(D_2\At D_2) - \tnorm{\mathcal{F}}} 
  \le \frac{\lambda_1(D_2\At D_2) + c(n,u)}{
  \lambda_n(D_2\At D_2) - c(n,u)}.
\end{equation}
$(D_2\At D_2)_{ii} = 1$ implies $\lambda_1(D_2\At D_2)\ge 1$,
which gives
\begin{equation}\notag
  \scond(\Atcomp) \le \scond(\At)
  \left( \frac{1 + \frac{c(n,u)}{\lambda_1(D_2\At D_2)}}{
  1 - \frac{c(n,u)}{\lambda_n(D_2\At D_2)}} \right)
  \le \scond(\At)
  \left( \frac{1 + c(n,u)}{1-c(n,u)\scond(\At)} \right),
\end{equation}
provided that $c(n,u) \scond(\At) < 1$. 
\end{proof}
The assumption that $\gamma_h$ satisfies Assumption~\ref{it.ass-size-gammah} in
Lemma~\ref{lem.bound-scAtcomp-scAt} is just for conciseness. The lemma is still valid without this assumption as long as
the $2u$ are replaced with 
$u + 4\gamma_h\alpha$ in the definition of $c(n,u)$ in~\eqref{def.cnu}.

Let us focus on $c(n,u)$ in~\eqref{def.cnu},  which is defined as an upper bound on 
$\tnorm{\mathcal{F}}$
in the proof of Lemma~\ref{lem.bound-scAtcomp-scAt}.
The factor of $n$ is pessimistic:
a simple numerical experiment that generates $500$ matrices in MATLAB with same size $n=10^2$, condition number $10^{13}$ and $5$ different singular values distributions, shows that $\tnorm{\mathcal{F}} \approx u$ in practice.
By assumption, $u < 0.1$ so $c(n,u) < 7nu$.
Now, let us go back to the bound~\eqref{bound-scAtcomp-scAt}.
If we enforce $7nu \scond(\At) < 1/2$
(Assumption~\ref{it.c(n,u)kappa<1/2}), then we have%
\begin{equation}\notag
  \scond(\Atcomp) \le
  \left( \frac{1+c(n,u)\scond(\At)}{1-c(n,u)\scond(\At)} \right)
  \scond(\At)
  \le
  \left( \frac{1+7nu\scond(\At)}{1-7nu\scond(\At)} \right)
  \scond(\At)
  \le 3\scond(\At).
\end{equation}
Substituting this inequality into \eqref{corollary-grand-error-bound} yields%
\begin{equation}\notag
  \ferrk \le 3 p_2 u \scond(\At) + 4nu\scond(\At) + p_1u
  = p_3\scond(\At)u + p_1u,
\end{equation}
where $p_3 := 3p_2 + 4n$,
and we have proved Theorem~\ref{thm.main}.

\subsection{Bounds on $\scond(\At)$}
\label{sec.disc-cond-numb}

\subsubsection{If $\At$ is scaled diagonally dominant} 
\label{sec.sdd-cond-numb}

A matrix $\At$ is \emph{$\theta$-scaled diagonally dominant}
($\theta$-s.d.d.) with respect to $\tnorm{\,\cdot\,}$
if $\tnorm{\wt D \At \wt D - I} \le \theta < 1$,
where $\wt D = \diag(\at_{ii}^{-1/2})$~\cite[sec.~2]{bade90}.

\begin{theorem}
[{\cite[Lem.~4.2, Prop.~4.3]{stwu02}}]
\label{thm.orth-bound}
Let $A \in \R\nn$.
If $\tnorm{A - I} \le \omega < 1$, then
$\kappa_{2}(A) \le (1+\omega)/(1-\omega)$. 
\end{theorem}

Consequently, if $\At$ is $\theta$-s.d.d., then 
\begin{equation}\label{eq.bound-on-kappa-sdd}
  \scond(\At) \le {1 + \theta\over 1 - \theta}.
\end{equation}
Let us now derive a condition on $\At$ such that
$\At$ is $\theta$-s.d.d.
We first write 
\begin{equation}\notag
  \tnorm{\Dt\At\Dt - I}=
  \tnorm{\Dt(\At - \Dt^{-2})\Dt}
  \le \tnorm{\Dt}^{2} \fnorm{\At - \Dt^{-2}}.
\end{equation}
Since $\Dt^{-2} = \diag(\at_{ii})$,
we can rewrite the last term as
$\tnorm{\Dt}^{2} \fnorm{\At - \diag(\at_{ii})}$
which is just $\off(\At)/\min_{i}(\at_{ii})$.
Therefore,
if the off-diagonals of $\At$ are sufficiently small compared to
the smallest diagonal entry of $\At$, such that
$\off(\At)/\min_{i}(\at_{ii}) = \theta < 1$,
then $\At$ is $\theta$-s.d.d. This implies, for example, that if $\off(\At)/\min_{i}(\at_{ii}) < 1/2$,
then $\scond(\At)$ is bounded by~$3$.

\subsubsection{The general case}
\label{sec.cdn-of-general-mtx}

The following result from Demmel and Veseli\'{c}~{\cite[Prop.~6.2]{deve92}
provides a bound on the scaled condition number.

\begin{theorem}\label{thm.hadamard-measure}
Let $d_{i}(A)$ and $\lambda_{i}(A)$ be
the $i$th largest diagonal entry and eigenvalue of a symmetric positive definite $A\in\R\nn$,
and let $D = \diag( a_{ii}^{-1/2} )$.
Then,%
\begin{equation}\notag
  \scond(A) = 
  \kappa_{2}(DAD) \le \frac{n}{\lambda_{n}(DAD)}
  \le n\, \eu \prod_{i=1}^{n} \frac{d_{i}(A)}{\lambda_{i}(A)},
\end{equation}
where $\eu$ is Euler's constant.
\end{theorem}
While the bound in Theorem~\ref{thm.hadamard-measure} is much less informative than that described in section~\ref{sec.sdd-cond-numb}, it demonstrates that it is still possible to have a modest value of $\scond(\At)$ even if $\tilde{A}$ is not $\theta$-s.d.d.~for any $\theta < 1$.

\section{Construction of the preconditioner}
\label{sec.altern-constr-prec}
We discuss several methods that 
exploit a lower precision $\ul > u$
in constructing a preconditioner $\Qt$ satisfying the following two properties.
\begin{enumerate}[label=\textup{(P\arabic*)}]
\item \label{it.preconditioner-P1}%
$\tnorm{\Qt\tp\Qt - I} \le p_{1}u<1$
\item \label{it.preconditioner-P2}%
$\off(\At) \le p_4\ul\fnorm{A}$.
\end{enumerate}

Property~\ref{it.preconditioner-P1} is required by Algorithm~\ref{alg.uh-prec-jacobi}, and together with Ostrowski's theorem \cite{ostrowski1959quantitative},
ensures that the \evals\ of $\At$ lie within a factor of $1 \pm p_1 u$ of those of $A$.
Property~\ref{it.preconditioner-P2}
benefits the convergence of the Jacobi algorithm
and is related to $\scond(\At)$ (see section \ref{sec.sdd-cond-numb}), which controls the relative forward error of the eigenvalues returned by Algorithm \ref{alg.uh-prec-jacobi}.

\subsection{Orthogonalizing an approximate \evec\ matrix}
\label{sec.method-A-orthogonal-evecmtx}
This approach is motivated by the fact that 
a low-precision spectral decomposition can approximate the one computed at working precision $u$. However, the eigenvector matrix is not orthogonal at precision $u$, so it needs orthogonalizing. 

\begin{algorithm2e}[h]
\caption{Constructing the preconditioner by orthogonalizing an
  \evec\ matrix computed at low precision $\ul$. }
\label{alg.orth-preconditioner}
\kwin{A symmetric matrix $A\in\R\nn$ and
  two precisions $\ul > u$. }
\kwout{A matrix $\Qt\in\R\nn$ satisfying~\ref{it.preconditioner-P1}
  and~\ref{it.preconditioner-P2}.}
\BlankLine
Compute a spectral decomposition $A = Q_{\ell}\Lambda_{\ell}Q_{\ell}\tp$
at precision $\ul$.\;
Orthogonalize $Q_{\ell}$ to $\Qt$ at precision $u$. \label{it.orth}
\end{algorithm2e}

Regarding the orthogonalization in line~\ref{it.orth},
Zhang and Bai~\cite{zhba22} use the modified Gram--Schmidt (MGS) method
and Zhou~\cite{zhou22-MPJA} uses the Newton--Schulz (NS) iteration~\cite[Eq.~8.20]{high08-FM}. We will also consider Householder QR factorization (HHQR).

\subsection{Using a low-precision tridiagonal reduction}
\label{sec.method-B-low-prec-trid}
We propose an alternative approach to Algorithm~\ref{alg.orth-preconditioner}
that does not require explicit reorthogonalization.
The key motivation is that for an arbitrary vector $v\in\R\n$,
the associated Householder reflector $H = I - 2vv\tp/\tnorm{v}^{2}$
is orthogonal to the precision of computation~\cite[Lem.~19.3]{high02-ASNA2},
regardless of the precision in which $v$ was initially computed or stored.

\begin{algorithm2e}[h]
\caption{Constructing the preconditioner using a low
  precision tridiagonalization.}
\label{alg.modified-tridiag}
\kwin{A symmetric matrix $A\in\R\nn$ and two precisions
  $\ul > u$. }
\kwout{A matrix $\Qt\in\R\nn$ satisfying~\ref{it.preconditioner-P1}
  and~\ref{it.preconditioner-P2}. }
\BlankLine
Tridiagonalize $A$ to $T_{\ell}$ using the Householder method at precision $\ul$. Store the Householder vectors in columns of $V \in \mathbb{R}^{n\times n-2}$.\;\label{li.methodB.tridiagonal}
Promote $V$ to precision $u$.\;
Compute a spectral decomposition $Q_{S} \Lambda_{\ell} Q_{S}\tp$ of $T_{\ell}$ at precision $u$.\;
Form $\wt Q$ by applying the Householder reflectors, constructed using columns of $V$, to $Q_S$ at precision $u$.\;\label{li.tridiag-4}
\end{algorithm2e}
In line~\ref{li.methodB.tridiagonal} of Algorithm~\ref{alg.modified-tridiag},
we perform the tridiagonal reduction at low precision,
since it is the most expensive part in terms of execution time of a symmetric eigensolver.

Note that most tridiagonal reduction routines, such as \texttt{SSYTRD} in LAPACK~\cite{abbb99-LUG3}, 
treat the Householder reflector as $H=I-\tau vv\tp$ with $\tau := 2/\tnorm{v}^2$.
Namely, in line~\ref{li.methodB.tridiagonal}, the tridiagonal reduction would give two outputs, 
$V$ and $\mathcal{T} \in \R^{n-2}$ where $\mathcal{T}_i$ is the corresponding constant for the $i$th Householder vector.
However, one should not use $\mathcal{T}$ in line~\ref{li.tridiag-4} directly since they are computed at low precision, which will cause $\Qt$ not to be orthogonal at working precision. 
Instead, one would need to recompute $\mathcal{T}$ using $V$ before applying them to $Q_S$.

\subsection{Proof of Properties \ref{it.preconditioner-P1} and \ref{it.preconditioner-P2}}
\label{sec.roundoff-methodA}

We need the following lemma on $\Qt$ computed by Algorithm~\ref{alg.orth-preconditioner}.
\begin{lemma}
\label{lem.bound-Qt-Ql}
Let $Q_{\ell}$ satisfy $\tnorm{Q_{\ell}\tp Q_{\ell} - I} \leq p_{5}u_\ell < 1/2$
and $\Qt$ be the orthogonalization of $Q_{\ell}$ using MGS, HHQR or NS at precision $u<\ul$. Then $\tnorm{Q_{\ell} - \Qt} \le p_{7} \ul$ and there exists $\Delta\Qt$ such that $\Qt+\Delta\Qt$ is orthogonal and $\tnorm{\Delta \Qt} \leq p_1 u$. 
\end{lemma}

\begin{proof}
The case of MGS was proved in~\cite{zhba22}. 
We start with the case of NS.
Let $\Qt$ be the computed unitary polar factor of $\Ql$ by applying NS. Since NS is a stable iteration when $\tnorm{\Ql}$ is safely bounded by $\sqrt{3}$~\cite[sec.~6.3]{nakatsukasa2012backward},
there exists a $\Delta \Qt$ such that $\Qt + \Delta\Qt = \Qpolar$ with $\tnorm{\Delta\Qt} \le p_1 u$, 
where $\Qpolar$ is the exact unitary polar factor of $\Ql + \Delta \Ql$ with 
$\tnorm{\Delta\Ql} \le p_6 u$~\cite{nakatsukasa2012backward}.
Using~\cite[Lem.~8.17]{high08-FM} and the bound on $\tnorm{\Ql\tp\!\Ql-I}$, we have%
\begin{equation}\label{eq.Ql-Qpolar}
  \tnorm{(\Ql+\Delta\Ql) - \Qpolar} \le \frac{\tnorm{(\Ql+\Delta\Ql)\tp(\Ql+\Delta\Ql) - I}}{
  1 + \sigma_{n}(\Ql+\Delta \Ql)} \le p_{5}\ul + O(u).
\end{equation}
Therefore, 
\begin{align}
    \label{eq.Ql-Qt-factorization} \tnorm{\Ql-\Qt} 
    & \le \tnorm{\Ql+\Delta\Ql-(\Qt+\Delta\Qt)} + \tnorm{\Delta\Ql}+\tnorm{\Delta\Qt} \\
    & \le p_5\ul+O(u) + p_6u + p_1u \le p_7\ul. \notag 
\end{align}

Now consider the case of HHQR.
Similarly to NS, there exists a $\Delta\Qt$ such that 
$\Qt+\Delta\Qt=\QQR$ with $\tnorm{\Delta\Qt} \le p_1 u$, 
where $\QQR$ is the exact orthogonal factor of 
the QR factorization of $\Ql+\Delta\Ql$ 
with $\tnorm{\Delta\Ql} \le p_6 u$~\cite[Eq.~19.13 \& Thm.~19.4]{high02-ASNA2}.
By \cite[Rmk~2.1]{sun95}, and~\eqref{eq.Ql-Qpolar}, we have
\begin{equation}
    \notag 
    \tnorm{\Ql+\Delta\Ql-\QQR} \le 5n^{1/2}\tnorm{\Ql+\Delta\Ql-\Qpolar} 
    \le 5n^{1/2}p_5\ul + O(u). 
\end{equation}
By expanding $\tnorm{\Ql-\Qt}$ in a similar way to~\eqref{eq.Ql-Qt-factorization}, we obtain $\tnorm{\Ql-\Qt} \le p_7\ul$. 
\end{proof}

Properties \ref{it.preconditioner-P1} and \ref{it.preconditioner-P2} are corollaries of the following backward error analysis.  
\begin{lemma}
\label{lem.property-preconditioner}
Let $A \in \R\nn$ be symmetric and let $\Qt$ be computed by Algorithm~\ref{alg.orth-preconditioner} or Algorithm~\ref{alg.modified-tridiag}. Then
\begin{equation}\label{eq.preconditioner-roundoff}
  A + \Delta A =
  (\Qt + \Delta \Qt) \Lambda_{\ell} (\Qt + \Delta \Qt)\tp,
\end{equation}
with $\Qt + \Delta \Qt$ orthogonal,
$\tnorm{\Delta\Qt} \le p_{8}u$ and
$\tnorm{\Delta A} \le p_{8} \ul \tnorm{A}$.
\end{lemma}

\begin{proof}
We divide the proof into two parts corresponding to the two algorithms. In the end, one can take $p_{8}$ as the maximum value of $p_{8}$ found for each algorithm.

Let us start with Algorithm~\ref{alg.orth-preconditioner}.
Since $Q_{\ell}\Lambda_{\ell}Q_{\ell}\tp$ is a spectral decomposition of $A$ computed at precision $u_{\ell}$, \cite[sec.~4.7.1]{abbb99-LUG3} gives
\begin{equation}\label{eq.roundoff-low-spec}
  A + \E_1 = (Q_{\ell} + \Delta Q_{\ell}) \Lambda_{\ell}(Q_{\ell} + \Delta Q_{\ell})\tp,
\end{equation}
with $Q_{\ell} + \Delta Q_{\ell}$ orthogonal,
$\tnorm{\E_1}\le p_{9} \ul \tnorm{A}$ and $\tnorm{\Delta Q_{\ell}} \le p_{9} \ul$.
Let $\Delta \Qt$ be as in Lemma \ref{lem.bound-Qt-Ql}. Extract the term $(\Qt+\Delta\Qt)\Lambda_\ell(\Qt+\Delta\Qt)\tp$
out of $A+\E_1$ to obtain
\begin{align}
  A + \E_1
  & = (\Qt + \Delta \Qt) \Lambda_{\ell} (\Qt + \Delta \Qt)\tp \notag\\
  & + (Q_{\ell} - \Qt + \Delta Q_{\ell} - \Delta \Qt) \Lambda_{\ell} 
    (Q_{\ell} + \Delta Q_{\ell})\tp \notag \\
  & + (\Qt + \Delta \Qt) \Lambda_{\ell}
    (Q_{\ell} - \Qt + \Delta Q_{\ell} - \Delta \Qt)\tp.
    \notag
\end{align}
Let $\E_2$ be the sum of the last two terms.  
Since $\Qt+\Delta\Qt$ and $\Ql+\Delta\Ql$ are orthogonal,%
\begin{equation}\label{eq.boundE2}
  \tnorm{\E_2}
  \le 2\tnorm{Q_{\ell} - \Qt + \Delta Q_{\ell} - \Delta \Qt}\tnorm{\Lambda_{\ell}},
\end{equation}
where
\begin{equation}\label{eq.error}
  \tnorm{Q_{\ell} - \Qt + \Delta Q_{\ell} - \Delta \Qt}
  \le \tnorm{Q_{\ell} - \Qt} + \tnorm{\Delta Q_{\ell}} + \tnorm{\Delta \Qt}
  \le p_{7}\ul+p_{9}\ul+p_1u,
\end{equation}
by Lemma \ref{lem.bound-Qt-Ql}. Letting $p_{10}\ul = p_{7}\ul + p_{9}\ul+p_1u$, and substituting~\eqref{eq.error} back into~\eqref{eq.boundE2}, we have 
$$
    \tnorm{\E_2} 
    \le 
    2p_{10}\ul\tnorm{\Lambda_\ell} =
    2p_{10}\ul\tnorm{A+\E_1} 
    \le 
     (2p_{10} + 2p_{10}p_{9}\ul)\ul\tnorm{A}.
$$
Define $p_{11}=2p_{10} + 2p_{10}p_{9}\ul$ and set $\Delta A = \E_1 - \E_2$. Then
\begin{equation}\notag
    A + \Delta A = A + (\E_1 - \E_2) = 
    (\Qt+\Delta \Qt)\Lambda_\ell(\Qt + \Delta \Qt)\tp,
\end{equation}
where $\tnorm{\Delta A} \le \tnorm{\E_1}+\tnorm{\E_2} \le (p_{9}+p_{11})\ul$ 
and $\tnorm{\Delta\Qt} \le p_1u$. Take $p_{8} = \max\{p_1,p_{9}+p_{11}\}$ to prove the lemma for Algorithm~\ref{alg.orth-preconditioner}.

Algorithm~\ref{alg.modified-tridiag}
starts with a tridiagonal reduction at precision $\ul$, which has the following backward error property~\cite[p.~297]{wilk65-AEP}:
\begin{equation}\label{eq.proof.alg3.tri}
  A + \E_3 = \wt Q_T  T_{\ell} \wt Q_T\tp,
\end{equation}   
with $\wt Q_T$ orthogonal 
(the exact product of the Householder reflectors used in the reduction)
and $\tnorm{\E_3} \le p_{12}\ul\tnorm{A}$~\cite[p.~297]{wilk65-AEP}.
In addition, the computed spectral decomposition of $T_{\ell}$ satisfies a relation similar to that in~\eqref{eq.roundoff-low-spec},
\begin{equation}\label{eq.proof.alg3.Tl}
  T_{\ell} + \Delta T_{\ell} = \wt Q_S \Lambda_{\ell} \wt Q_S\tp,
\end{equation}
with $\wt Q_S=Q_S+\Delta Q_S$ orthogonal,
$\tnorm{\Delta Q_S}\le p_{13}u$, and $\tnorm{\Delta T_{\ell}}\le p_{13}u\tnorm{T_{\ell}}$.
Substituting~\eqref{eq.proof.alg3.Tl}  back into~\eqref{eq.proof.alg3.tri} gives
\begin{equation}\label{roundoff-1}
  A+\Delta A = G \Lambda_{\ell} G\tp,
\end{equation}
where $G=\wt Q_T\wt Q_S$ is orthogonal and
$\Delta A = \E_3+ \wt Q_T \Delta T_{\ell} \wt Q_T\tp$
with 
$$
\tnorm{\Delta A}\le p_{12} \ul\tnorm{A} + p_{13} u\tnorm{T_{\ell}}=
p_{14} \ul\tnorm{A}.
$$

Finally, let $\Qt$ be the matrix obtained by applying, at precision $u$, the Householder reflectors used in the reduction to tridiagonal form to $Q_S$. Then it follows from \cite[Lem.~19.3]{high02-ASNA2} that
\begin{equation}\label{eq.proof.alg3.mmp}
  \Qt + \E_4 = \wt Q_T Q_S, \qquad \tnorm{\E_4} \le p_{14} u \tnorm{Q_S}
  = p_{15}u.
\end{equation}
But $Q_S = \Qt_S-\Delta Q_S$ with $\tnorm{\Delta Q_S}\le p_{13}u$ so that if we let $\Delta \Qt=\E_4+\wt Q_T\Delta Q_S$ then
$\Qt+\Delta \Qt = \wt Q_T \wt Q_S=G$ is orthogonal  with 
$$
\tnorm{\Delta\Qt}\le p_{13}u+p_{15}u=p_{16}u.
$$
Taking $p_{8} = \max\{p_{14}, p_{16}\}$ completes the proof. 
\end{proof}

\begin{corollary}
    Let $A \in \R\nn$ be symmetric, and let $\Qt$ be computed by Algorithm~\ref{alg.orth-preconditioner} or Algorithm~\ref{alg.modified-tridiag}. Then $\Qt$ satisfies Properties \ref{it.preconditioner-P1} and \ref{it.preconditioner-P2}.
\begin{proof}
    We start with Property \ref{it.preconditioner-P1}. By Lemma \ref{lem.property-preconditioner}, $\Qt + \Delta \Qt$ is orthogonal for some matrix $\Delta\Qt$ satisfying $\tnorm{\Delta\Qt} \leq p_{8} u$. This implies
    \begin{align*}
        I - \Qt^T \Qt &= (\Qt+\Delta\Qt)^T (\Qt+\Delta\Qt) - \Qt^T\Qt \\
        &= \Delta\Qt^T (\Qt+\Delta\Qt) + (\Qt+\Delta\Qt)^T \Delta \Qt - \Delta\Qt^T\Delta\Qt.
    \end{align*}
    Taking the 2-norm results in
    \begin{equation*}
        \tnorm{\Qt^T \Qt - I} \leq 2p_{8}u + (p_{8} u)^2 \leq p_1u.
    \end{equation*}

    Now we prove Property \ref{it.preconditioner-P2}. Observe that $\At-\diag(\At)$ has zero diagonal entries and the same off-diagonal entries as $\Qt\tp\!A\Qt-\Lambda_\ell$. Therefore,
\begin{equation}\notag
  \off(\At) = \fnorm{\At - \diag(\At)}
  \le \fnorm{\Qt\tp A \Qt - \Lambda_{\ell}}.
\end{equation}
By Lemma~\ref{lem.property-preconditioner},
\begin{equation}\notag
  (\Qt + \Delta \Qt)\tp A (\Qt + \Delta \Qt) +
  (\Qt + \Delta \Qt)\tp \Delta A (\Qt + \Delta \Qt) = \Lambda_{\ell},
\end{equation}
where $\tnorm{\Delta A} \leq p_{8} u_\ell \|A\|_2$. This implies 
\begin{align*}
    \Lambda_\ell-\Qt\tp\!A\Qt &= \Delta \Qt ^T A (\Qt+\Delta\Qt) + (\Qt+\Delta\Qt)^T A\Delta\Qt - \Delta\Qt^T A \Delta\Qt \\
    & \qquad + (\Qt+\Delta\Qt)\tp \Delta A (\Qt+\Delta\Qt).
\end{align*}
Taking the $2$-norm results in
\begin{align*}
\tnorm{\Qt\tp A \Qt - \Lambda_{\ell}}
& \le 2 \tnorm{\Delta \Qt}\tnorm{A} + \tnorm{\Delta\Qt}^2 \tnorm{A} + \tnorm{\Delta A} \\
& \le 2p_{8} u\tnorm{A} + (p_{8}u)^2 \tnorm{A} + p_{8} \ul \tnorm{A} \\
& \le p_{17} \ul \tnorm{A}.
\end{align*}
Since $\fnorm{A} \le n^{1/2}\tnorm{A}$, we have 
\begin{equation}\notag
\fnorm{\Qt\tp A \Qt - \Lambda_{\ell}} 
\le n^{1/2}\tnorm{\Qt\tp A \Qt - \Lambda_{\ell}}
\le n^{1/2}p_{17}\ul\tnorm{A}\le n^{1/2}p_{17}\ul\fnorm{A}.
\end{equation}
Letting $p_{4} = n^{1/2}p_{17}$, we have shown that the preconditioners constructed using Algorithms~\ref{alg.orth-preconditioner} and~\ref{alg.modified-tridiag} satisfy property~\ref{it.preconditioner-P2}.
\end{proof}
\end{corollary}

\begin{remark}\label{rmk.choice of ul}
In section~\ref{sec.disc-cond-numb} we showed that if ${\off(\At)}/{\min_i(\at_{ii})} = \theta < 1$, then $\scond(\At) \le (1+\theta)/(1-\theta)$.
Combining Property \ref{it.preconditioner-P2} with Theorem~\ref{thm.majorization} and Ostrowski's Theorem \cite{ostrowski1959quantitative}, we obtain
\begin{equation}\notag
  \frac{\off(\At)}{\min_{i}(\at_{ii})} \le
  \frac{p_{4}\ul\fnorm{A}}{\lambda_{n}(\At)}
  \le p_{18}\ul\kappa_{2}(A),
\end{equation}
where $p_{18} = n^{1/2} p_4 / (1 - p_1u)$. Therefore, we can bound $\scond(\At)$ if $\kappa_2(A) < 1/(p_{18}\ul)$.
\end{remark}

\section{Numerical Experiments}
\label{sec.numerical-examples}

We performed numerical experiments to assess the quality of the preconditioners, the accuracy of the computed eigenvalues, and the magnitude of the scaled condition number $\scond(\At)$ (which is relevant to the forward error bound of our main theorem) using MATLAB R2024b. In addition, we ran a performance test using Julia (version 1.11.5).
We fixed $\ul$, $u$, and $\uh$ as IEEE single, double, and quadruple precisions, respectively.
Hence, $\ul = 2^{-24}$, $u = 2^{-53}$, and $\uh = 2^{-113}$.
To avoid overflow, we follow the scaling techniques in LAPACK routines for symmetric eigenvalue problems (see the code for \texttt{DSYEV}, and also \cite{hpz19}). We output a warning if the underflow is detected after demoting from $\Athcomp$ to $\Atcomp$.

To make our results reproducible, we used the commands \inline|rng(1)| (MATLAB) and \inline|Random.seed!(1)| (Julia) to control the random number generator.
We conducted our experiments using a 2023 MacBook Pro featuring an M3 Pro processor and 32 GiB of RAM. 
Quadruple precision was implemented using the MATLAB Multiprecision Computing Toolbox~\cite{mct2023} and the Julia package Quadmath (version 0.5.13).

Some of the experiments involve random dense \spd\ matrices. 
We generated such a matrix $A\in\R\nn$ using the MATLAB test matrices collection \inline|A = gallery('randsvd', N, -KAPPA, MODE);| where $\texttt{N} = n$, $\texttt{KAPPA} = \kappa_{2}(A)$, and $\texttt{MODE}=\text{MODE} \in \{1,2,\dots,5\}$, which specifies the singular value distributions of $A$. We will not go into the details of these different distributions; for more detailed information see \url{https://mathworks.com/help/releases/R2024b/matlab/ref/gallery.html}.
Code to reproduce the results of this section can be found at \url{https://github.com/zhengbo0503/Code_htwz25}.

\subsection{The size of the off-quantity}
\label{sec.numexp.sizeofoffA}
We investigated the size of $\off(\At)$ relative to $\|A\|_F$. 
The upper bound we derived on $\off(\At)$ is independent of the condition number of $A$, so the test matrices in this experiment have varying $n$ and MODE, but fixed condition number $\kappa_{2}(A) = 10^{6}$.
The results reported in Figure~\ref{fig.offAt} are in alignment with the theory developed in section~\ref{sec.altern-constr-prec}. With the exception of MODE = 2, we see a clear $n^{1/2}$ growth in the off-quantity. We observe that for MODE = 1, 2 and 4, there is no clear difference for the off-quantity between the different algorithms.
However, for MODE = 3 and MODE = 5, 
the off-quantity for Algorithm~\ref{alg.modified-tridiag} is smaller 
than that corresponding to Algorithm~\ref{alg.orth-preconditioner}.
This phenomenon remains for different RNG seeds,
and we currently do not have an explanation for it.

\begin{figure}[h!]
\centering \footnotesize
\subfloat[MODE = 1.]{\input{figs/offAt_mode1.tex}} \qquad\qquad
\subfloat[MODE = 2.]{\input{figs/offAt_mode2.tex}}

\subfloat[MODE = 3.]{\input{figs/offAt_mode3.tex}} \qquad\qquad
\subfloat[MODE = 4.]{\input{figs/offAt_mode4.tex}}

\hspace{18pt}
\subfloat[MODE = 5.]{\input{figs/offAt_mode5.tex}}
\hspace{1pt}
\subfloat{\input{figs/offAt_legend.tex}}
\caption{Size of $\off(\At)/\fnorm{A}$ against the matrix size $n$ of $\At=\wt Q\tp A\wt Q$ with $\wt Q$ computed by four different algorithms specified in the legend. 
The subcaptions state the MODE used in \inline|randsvd|.
In all cases $\kappa_{2}(A) = 10^{6}$.}
\label{fig.offAt}
\end{figure}

\subsection{Relative forward error}
\label{sec.relat-forw-error}
We investigated the relative forward error bound of section~\ref{sec.main-result-proof}.
In the first experiment, we varied $\kappa_{2}(A)$,
while in the second, we varied $n$. We compared four different algorithms:
\begin{enumerate}[label=\upshape(\arabic*)]
\item {\bf Jacobi:} 
Coded according to~\cite[Alg.~3.1]{deve92} with stopping tolerance $\tol$ set to $n^{1/2}u$.
This choice is consistent with the LAPACK routine \texttt{DGESVJ}.
\item {\bf MP2Jacobi:}
Algorithm~\ref{alg.uh-prec-jacobi} with $\uh$ set to $u$,
together with Algorithm~\ref{alg.orth-preconditioner}
with HHQR to construct the preconditioner $\Qt$.
The term ``MP2'' emphasizes that this is a mixed-precision algorithm
with two precisions $\ul$ and $u$.
\item {\bf MP3Jacobi:}
Algorithm~\ref{alg.uh-prec-jacobi} together with Algorithm~\ref{alg.orth-preconditioner} with HHQR to construct the preconditioner $\Qt$. The term ``MP3'' emphasizes that this is a mixed-precision algorithm with three precisions $\ul$, $u$, and $\uh$.
\item {\bf MATLAB:}
The MATLAB function \texttt{eig} which performs a tridiagonal reduction before applying a tridiagonal eigensolver~\cite[sec.~9.8.2]{hihi16-MG3}.
\end{enumerate}

We use Algorithm~\ref{alg.orth-preconditioner} with HHQR to construct the preconditioner in both MP2Jacobi and MP3Jacobi since HHQR has already been implemented as \texttt{qr} in MATLAB. Other constructions such as Algorithm~\ref{alg.orth-preconditioner} with MGS and NS, and Algorithm~\ref{alg.modified-tridiag} will produce similar results in this section. 

For each experiment, we plotted the error bound 
\begin{equation}
    \ferrk \le 7n\scond(\At) u,
\end{equation}
which results from absorbing the second term on the right-hand side of~\eqref{eq:ferrk-main} into the first, and setting $p_2 = 7n$. 
We approximated the exact $\At$ by computing the product $\Qt\tp\!A\Qt$  at octuple precision to compute this bound.

\subsubsection{Varying $\kappa_2(A)$}
\label{sec.verify-bound-with-varying-2norm}
Figure~\ref{fig.maxferrk-varykappa} displays
the maximal relative forward error $\max_{k} \ferrk$ and the bound $7n\scond(\At)u$ for matrices $A \in \R^{100\times 100}$ with condition numbers between $10$ and $10^{16}$.
The plots show that the algorithm MP3Jacobi produced the smallest maximum \rfe\ among all four algorithms.
Even for matrices with a condition number of $10^{16}$,
it still secured about 8 digits of accuracy.
For well-conditioned matrices,
all four algorithms attained a similar level of accuracy.

\begin{figure}[h!] 
\centering\footnotesize
\subfloat[MODE = 1.]{\input{figs/varykappa_mode1.tex}} \qquad\qquad
\subfloat[MODE = 2.]{\input{figs/varykappa_mode2.tex}}

\subfloat[MODE = 3.]{\input{figs/varykappa_mode3.tex}} \qquad\qquad
\subfloat[MODE = 4.]{\input{figs/varykappa_mode4.tex}}

\hspace{14pt}\subfloat[MODE = 5.]{\input{figs/varykappa_mode5.tex}}
\subfloat{\input{figs/varykappa_legend.tex}}

\caption{Maximum \rfe, $\max_{k}\ferrk$, against $\kappa_{2}(A)$
  for four different algorithms and 
  $A\in\R^{100\times 100}$ generated by \inline|randsvd| with different values of MODE.}
\label{fig.maxferrk-varykappa}
\end{figure}

\subsubsection{Varying the matrix size}
\label{sec.verify-bound-varying}
Figure~\ref{fig.maxferrk-varydim} displays 
$\max_{k} \ferrk$ as well as the bound $7n\scond(\At)u$ for $A \in \R^{n\times n}$ with $10\le n\le 10^{3}$ and $\kappa_2(A)=10^{10}$.
Just as in section \ref{sec.verify-bound-with-varying-2norm}, the algorithm MP3Jacobi produced the smallest maximum \rfe\ among all four algorithms.
For test matrices with $\text{MODE}\in\{3,5\}$, the bound $7n\scond(\At)u$ is pessimistic. This phenomenon also appears in Figure~\ref{fig.maxferrk-varykappa}, and we will see the cause in section~\ref{sec.reduct-cond-numb}.

\begin{figure}[h!]
\centering\footnotesize
\subfloat[MODE=1.]{\input{figs/varydim_mode1.tex}}\qquad\qquad
\subfloat[MODE=2.]{\input{figs/varydim_mode2.tex}}

\subfloat[MODE=3.]{\input{figs/varydim_mode3.tex}}\qquad\qquad
\subfloat[MODE=4.]{\input{figs/varydim_mode4.tex}}

\hspace{15pt}
\subfloat[MODE=5.]{\input{figs/varydim_mode5.tex}}
\hspace{1pt}
\subfloat{\input{figs/varydim_legend.tex}}

\caption{Maximum relative forward error,
  $\max_{k}\ferrk$, for four different algorithms against the  size $n$ of matrices $A$ generated by \inline{randsvd} with different values of MODE and $\kappa_{2}(A) = 10^{10}$.}
\label{fig.maxferrk-varydim}
\end{figure}

\subsubsection{Special test matrices}
\label{sec.evals-special-mtx}
We considered the following test matrices.
\begin{enumerate}[label=\upshape(\arabic*),labelwidth=!]
\item 
\inline|gallery('randsvd', 100, -1e8, 3)|: 
random symmetric positive definite matrix $A\in\R^{100\times 100}$
with $\kappa_2(A)=10^8$ and geometrically distributed singular values. 

\item 
\inline|anymatrix('nessie/whiskycorr')|: 
symmetric positive definite correlation matrix $A\in\R^{86\times 86}$ with $\kappa_2(A) = 4.9\times 10^{19}$ from the Anymatrix Toolbox~\cite{himi21}.

\item \label{it.special-mtx-cov}
Covariance matrix of a data matrix from~\cite{superconductivty_data_464} which describes several superconductors properties. This matrix is $82\times 82$ and has condition number   $2.1\times 10^{12}$. 
\item
Gram matrix of \inline|gallery('lauchli', 500, 1e-3)|.
The matrix is $500\times 500$ with condition number $5\times 10^8$.

\end{enumerate} 

\begin{figure}[h!]
\centering\footnotesize
\subfloat[\inline|gallery('randsvd', 100, -1e8, 3)|\label{subfig:random}]{\input{figs/rand_matrix.tex}}

\subfloat[\inline|anymatrix('nessie/whiskycorr')|\label{subfig:anymatrix}]{\input{figs/whisky.tex}}

\subfloat[Covariance matrix\label{subfig:covariance matrix}]{\input{figs/covariance.tex}}

\subfloat[Gram matrix of \inline|gallery('lauchli', 500, 1e-3)|\label{subfig.lauchli}]{\input{figs/nos7.tex}}

\input{figs/special_legend.tex}
\caption{Relative forward error for the eigenvalues ordered from largest to smallest of the test matrices in section~\ref{sec.evals-special-mtx}.
Subfigure (d)
only displays $\ferrk$ for every tenth eigenvalue, and the data missing for MP3Jacobi and MATLAB is because the error was smaller than $10^{-16}$.}
\label{fig.special-mtx}
\end{figure}

The first three plots in Figure~\ref{fig.special-mtx} display $\ferrk$ for each eigenvalue and the fourth only shows $\ferrk$ for every tenth eigenvalue to avoid clutter in the plot.
For these test matrices, the plots show that 
only MP3Jacobi consistently computes all eigenvalues with high relative accuracy. The relative forward errors for the other algorithms can behave differently, depending on the matrix in question. Figure~\ref{subfig:random} shows the ``typical" behavior: all eigensolvers compute the largest eigenvalues accurately, but MP2Jacobi, Jacobi, and MATLAB gradually lose accuracy for smaller eigenvalues.  
Figure~\ref{subfig:anymatrix} shows the loss of accuracy for the small eigenvalues computed by 
MP2Jacobi, Jacobi, and MATLAB on an extremely ill-conditioned matrix. Figure~\ref{subfig:covariance matrix} shows that MP2Jacobi and Jacobi can maintain a high relative accuracy if the ill-conditioning is solely due to poor scaling ($\kappa_2(A) \approx 10^{12}$, but $\scond(A) \approx 10^5$).
In Figure~\ref{subfig.lauchli} we observe that MATLAB computed extremely accurate eigenvalues, but less accurate than those computed by MP3Jacobi. However, MATLAB lost accuracy for the last few eigenvalues. MP2Jacobi and Jacobi were uniformly much less accurate.


\subsection{The scaled condition number of the preconditioned matrix}
\label{sec.reduct-cond-numb}
We investigated how much the scaled condition number reduces after applying the preconditioners from section~\ref{sec.altern-constr-prec}. We generated $A\in\R^{100\times 100}$ with condition numbers between $10$ and $10^{12}$, and five different singular value distributions determined by \texttt{randsvd}'s MODE. We report the ratios $\scond(A)/\k(A)$ and $\scond(\At)/\k(A)$ in Figure~\ref{fig.reduct_cond}, where $\At$ was computed in octuple precision.

\begin{figure}[h!]
\centering\footnotesize 
\subfloat[MODE=1.]{\input{figs/reducekappa_mode1.tex}}\qquad\qquad
\subfloat[MODE=2.]{\input{figs/reducekappa_mode2.tex}}

\subfloat[MODE=3.]{\input{figs/reducekappa_mode3.tex}}\qquad\qquad
\subfloat[MODE=4.]{\input{figs/reducekappa_mode4.tex}}

\qquad
\subfloat[MODE=5.]{\input{figs/reducekappa_mode5.tex}}\quad
\subfloat{\input{figs/reducekappa_legend.tex}}

\caption{Ratio of the scaled condition number of $A$ and $\At$, and $\kappa_{2}(A)$. The subcaptions state the singular value distributions for the test matrices $A$ which have fixed size $n=100$. }
\label{fig.reduct_cond}
\end{figure}

In Figure~\ref{fig.reduct_cond} we observe that $\scond(A)$ is of a similar magnitude to $\k(A)$ for these random matrices, while $\scond(\At)$ can be significantly smaller than $\k(A)$ when $A$ is even only moderately ill-conditioned. However, this effect is less pronounced for MODE = $3$ and $5$, appearing to stall at a ratio of $10^7$. This phenomenon could be related to the behavior of the forward error for these MODEs, but it is not clear to us at present.

In addition, we picked some ill-conditioned \spd\ test matrices and performed the same experiment. Table~\ref{table.anymatrix-test-reduct-kappa} shows our findings.

\begin{table}[h!]
\centering\footnotesize
\caption{Size of the scaled condition numbers
   $\scond(A)$ and $\scond(\At)$.}
\label{table.anymatrix-test-reduct-kappa}

\begin{tabular}{ccccc} \toprule Matrix Type & $n$ & $\kappa_{2}(A)$ & $\scond(A)$ & $\scond(\At)$\\ 
\midrule 
\inline|anymatrix('matlab/hilb')| & 20 & 1.08e+18 & 3.56e+18 & 3.66e+09 \\ \inline|anymatrix('matlab/invhilb')| & 20 & 2.14e+25 & 5.83e+17 & 3.26e+10 \\ \inline|anymatrix('matlab/pascal')| & 15 & 2.84e+15 & 5.83e+12 & 1.55e+04 \\ \inline|anymatrix('nessie/traincorr')| & 25 & 1.84e+18 & 9.97e+19 & 2.42e+05 \\
\inline|anymatrix('nessie/whiskycorr')| & 86 & 4.90e+19 & 1.12e+19 & 2.97e+06 \\ 
Covariance matrix (Section~\ref{sec.evals-special-mtx}) & 82 & 2.13e+12 & 6.46e+05 & 9.45e+04 \\ 
Gram matrix (Section~\ref{sec.evals-special-mtx}) & 500 & 5.00e+08 & 5.00e+08 & 1.09e+00 \\
\bottomrule \end{tabular}
\end{table}

\subsection{Timing with Julia}
\label{subsec.timing-test}
We present the results of timing Julia implementations of Jacobi, MP2Jacobi, and MP3Jacobi, available in the registered Julia package, JacobiEigen.jl.
Julia functions can be compiled into efficient implementations, which provide insight into how our algorithms would perform if implemented in an optimized library. 

\begin{figure}
    \centering
    \includegraphics[width=0.85\linewidth]{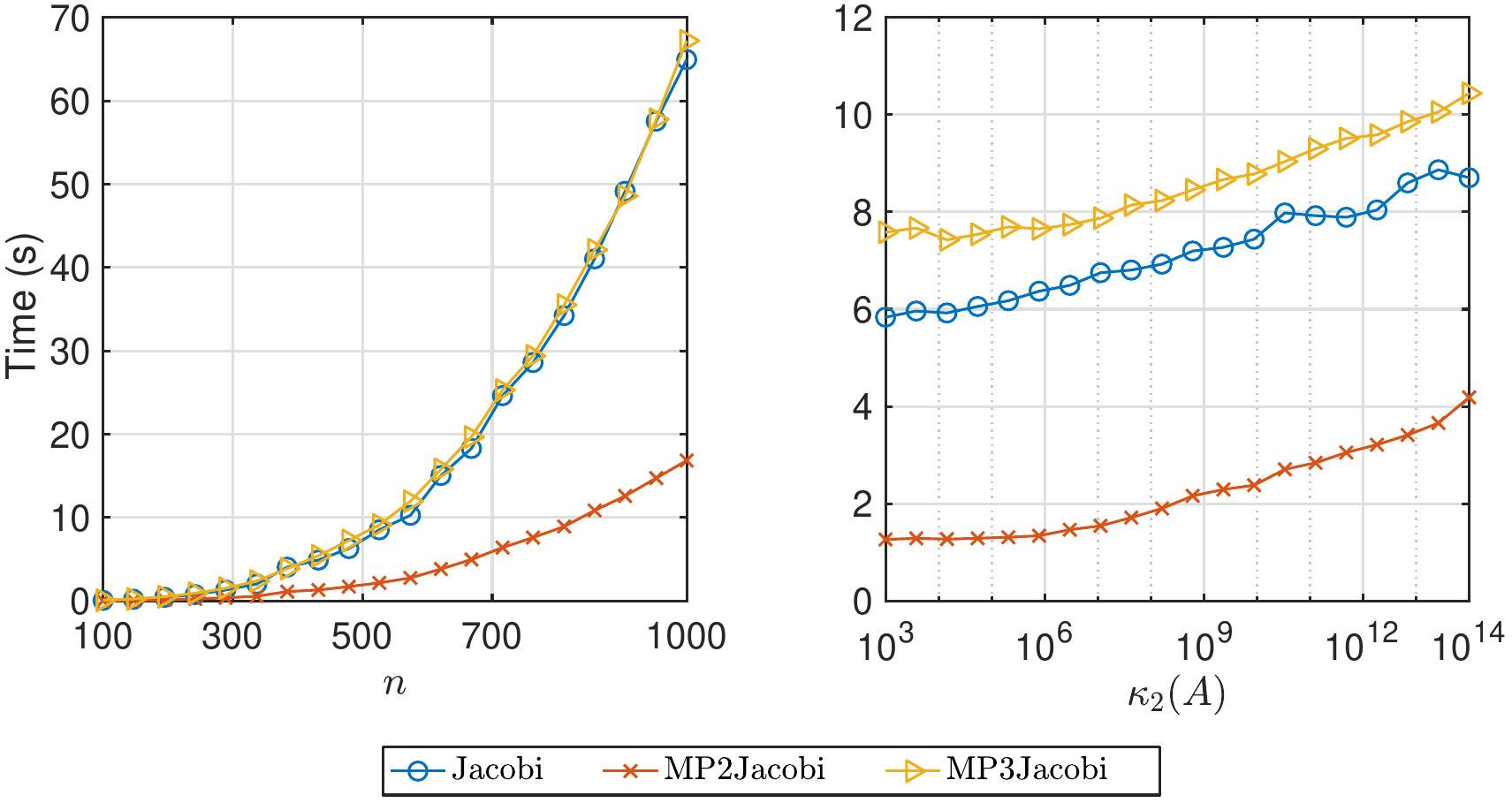}
    \caption{Execution time (in seconds) for performing Jacobi, MP2Jacobi, and MP3Jacobi. The test matrices on the left vary in size from $n=10^2$ to $n=10^3$ with a fixed condition number of $10^8$, while those on the right vary the 2-norm condition number from $10^3$ to $10^{14}$ with a fixed size of $n=500$.}
    \label{fig:timing-union}
\end{figure}

We constructed two sets of test matrices: one consisting of 20 matrices whose sizes vary from $n=10^2$ to $n=10^3$ with a fixed condition number of $10^8$, and another consisting of 20 matrices of size $n=500$ whose condition numbers range from $10^3$ to $10^{14}$. Figure~\ref{fig:timing-union} presents the runtimes of the Jacobi, MP2Jacobi, and MP3Jacobi algorithms on these matrices.

MP2Jacobi consistently outperforms both Jacobi and MP3Jacobi. For moderate matrix sizes, say $n \lesssim 500$, MP3Jacobi incurs roughly a $20\%$ overhead compared to Jacobi, but this relative overhead decreases as $n$ grows. As $n$ increases, the timings of Jacobi and MP3Jacobi become essentially indistinguishable, demonstrating that MP3Jacobi delivers improved accuracy with little extra cost. On the other hand, increasing the condition number produces a similar impact on the execution time of all three algorithms. This trend arises because, for ill-conditioned matrices, the Jacobi method tends to perform more sweeps to achieve convergence.

Figure \ref{fig.timing_all} isolates the timings of applying the preconditioner to $A$, and applying the Jacobi algorithm to the preconditioned matrix. As expected, it is only the former that severely affects the difference in timing. Indeed, the number of sweeps in the Jacobi algorithm depends only on $\off(\Atcomp)$ and the eigenvalue gap~\cite{hari91}, which is not significantly affected by the precision at which the preconditioner is applied (see Figure \ref{fig.timing_all} (b)).

Figure~\ref{fig.timing_all} (a) shows that applying the preconditioner at
quadruple precision was more than $1000$ times slower than at double
precision, which is an expected behavior due to quadruple precision not
being widely supported by hardware~\cite{knn23}. This is a facet of our
implementation --- in theory, it could be only about $7$ times slower, and
in practice, it is usually less than $50$ times slower, see \cite{babo15} \cite{bail05} and \cite{divi06}. 
As a result, 
Figure~\ref{fig.timing_all} (c) shows that our implementation takes roughly four times as long as MP2Jacobi because applying the preconditioner in quadruple precision dominates the runtime---a step that is negligible in MP2Jacobi.
Moreover, the red dashed line shows that if forming $\Atcomp$ at quadruple precision is only $100$ times slower compared to that at double precision, then MP3Jacobi (Potential) would take about the same time as MP2Jacobi. 
This implies that MP3Jacobi algorithm can be both fast and accurate if software or hardware support for quadruple-precision matrix-matrix multiplication improves.




\begin{figure}
    \centering
    \includegraphics[width=0.85\linewidth]{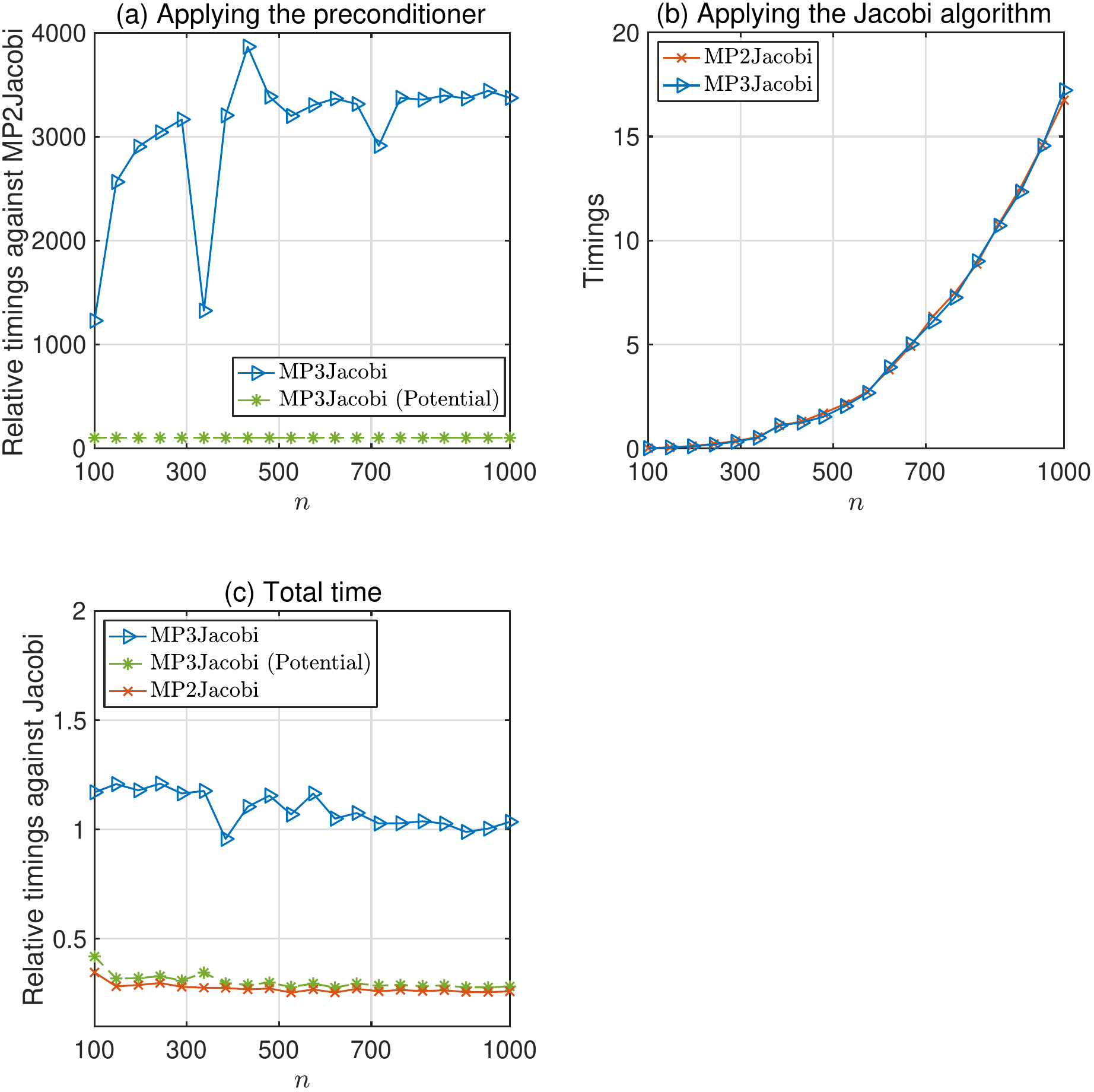}
    \caption{Timing results for MP2Jacobi, MP3Jacobi, and a theoretical variant of MP3Jacobi, MP3Jacobi (Potential), which represents a scenario where the preconditioner application stage in MP3Jacobi takes only $100$ units more time than in MP2Jacobi.}
    \label{fig.timing_all}
\end{figure}


\section{Conclusions and further remarks}
\label{sec.conclusion}
Demmel and Veseli\'{c} proved that with an appropriate stopping criterion, the Jacobi algorithm computes highly accurate eigenvalues of a symmetric positive definite matrix that is not properly scaled~\cite[Alg.~3.1]{deve92}, e.g.~a graded matrix. 
In this paper, we proved that the eigenvalues of a wider class of matrices are computed with high accuracy by adding an appropriate preprocessing step, which consists of applying a preconditioner at a higher precision $u_h$.

In addition, we collated several methods, one of which is new, for computing a preconditioner $\Qt$, that takes advantage of a lower precision $u_\ell$. We proved bounds on $\|\Qt^T \Qt - I\|_2$ and $\off(\At)$, relevant to $E_3$ from the proof of the forward error bound and convergence of the Jacobi algorithm~\cite{hari91}, respectively.

Our analysis suggests that a practical implementation of Algorithm~\ref{alg.uh-prec-jacobi} could be developed that dynamically chooses the lower and higher precisions. Regarding the higher precision, Assumption \ref{it.ass-size-gammah} suggests the choice $u_h = c_1 u \kappa_2(A)^{-1}$ for some hyperparameter $c_1$. Regarding the lower precision, by Remark \ref{rmk.choice of ul},
\begin{equation}
    u_\ell < \frac{1}{p_{18} \kappa_2(A)} \implies \frac{\off(\At)}{\min_i \tilde{a}_{ii}} = \theta < 1 \implies \scond(\At) < \frac{1+\theta}{1-\theta}.
\end{equation}
This suggests the choice $u_\ell = c_2 \kappa_2(A)^{-1}$ for a hyperparameter $c_2$. 
Since computing the exact condition number is costly, we recommend using
an estimate of $\kappa_2(A)$, which may be analytically or algorithmically available \cite{higham2000block, dixon1983estimating}. 

One could consider applying $\Qt$ to $A$ in factored form in line 1 of Algorithm~\ref{alg.uh-prec-jacobi}, rather than forming the matrix $\Qt$ explicitly and using matrix-matrix multiplications. 
We do not expect the error bounds to differ greatly, but we have not investigated this.

In a forthcoming paper, we will investigate a one-sided version of the algorithm to compute highly accurate singular values.


\section*{Acknowledgements}
We are grateful to Zlatko Drma\v{c} and Xiaobo Liu for fruitful discussions during the research for this paper, and to the anonymous reviewers for their insightful comments. We thank Thomas House for his advice on covariance matrices.
The work of the first author was supported by the Royal Society, the second author was supported by Engineering and Physical Sciences Research Council grant EP/W018101/1, and the fourth author was supported by  the University of Manchester Research
Scholar Award. 

\bibliographystyle{siamplain}
\bibliography{ref} 

\end{document}

%% file: figs/offAt_mode1.tex
%
\definecolor{mycolor1}{rgb}{0.00000,0.44700,0.74100}%
\definecolor{mycolor2}{rgb}{0.85000,0.32500,0.09800}%
\definecolor{mycolor3}{rgb}{0.92900,0.69400,0.12500}%
\definecolor{mycolor4}{rgb}{0.49400,0.18400,0.55600}%
\pgfplotstableread[col sep=comma]{figs/offAt_mode1.csv}{\offAtmodeI}

\begin{tikzpicture}[trim axis left,trim axis right]
  \begin{axis}[%
    xmode=log,
    xmin=10,
    xmax=1000,
    xminorticks=true,
    ymode=log,
    ymin=5e-08,
    ymax=1e-05,
    yminorticks=true,
    xlabel={$n$},
    ylabel={$\off(\wt{A})/\fnorm{A}$},
    ]
    \addplot[color=mycolor1, only marks,
    mark size=3.0pt, mark=triangle, mark options={solid, mycolor1}]
    table [x=n, y=qr] {\offAtmodeI};
    \addplot [color=mycolor2, only marks,
    mark size=3.0pt, mark=diamond, mark options={solid, mycolor2}]
    table [x=n, y=mgs] {\offAtmodeI};
    \addplot [color=mycolor3, only marks,
    mark size=3.0pt, mark=x, mark options={solid, mycolor3}]
    table [x=n, y=ns] {\offAtmodeI};
    \addplot [color=mycolor4, only marks,
    mark size=3.0pt, mark=star, mark options={solid, mycolor4}]
    table [x=n, y=tri] {\offAtmodeI};
    \addplot [color=black, dashed, line width=2.0pt]
    table [x=n, y=bound] {\offAtmodeI};
  \end{axis}
\end{tikzpicture}%


%% file: figs/offAt_mode2.tex
%
\definecolor{mycolor1}{rgb}{0.00000,0.44700,0.74100}%
\definecolor{mycolor2}{rgb}{0.85000,0.32500,0.09800}%
\definecolor{mycolor3}{rgb}{0.92900,0.69400,0.12500}%
\definecolor{mycolor4}{rgb}{0.49400,0.18400,0.55600}%

\pgfplotstableread[col sep=comma]{figs/offAt_mode2.csv}{\offAtmodeII}
\begin{tikzpicture}[trim axis left, trim axis right]
  \begin{axis}[%
    xmode=log,
    xmin=10,
    xmax=1000,
    xminorticks=true,
    ymode=log,
    ymin=1e-09,
    ymax=1e-05,
    yminorticks=true,
    xlabel={$n$},
    ]
    \addplot[color=mycolor1, only marks,
    mark size=3.0pt, mark=triangle, mark options={solid, mycolor1}]
    table [x=n, y=qr] {\offAtmodeII};
    \addplot [color=mycolor2, only marks,
    mark size=3.0pt, mark=diamond, mark options={solid, mycolor2}]
    table [x=n, y=mgs] {\offAtmodeII};
    \addplot [color=mycolor3, only marks,
    mark size=3.0pt, mark=x, mark options={solid, mycolor3}]
    table [x=n, y=ns] {\offAtmodeII};
    \addplot [color=mycolor4, only marks,
    mark size=3.0pt, mark=star, mark options={solid, mycolor4}]
    table [x=n, y=tri] {\offAtmodeII};
    \addplot [color=black, dashed, line width=2.0pt]
    table [x=n, y=bound] {\offAtmodeII};
  \end{axis}
\end{tikzpicture}%


%% file: figs/offAt_mode3.tex
%
\definecolor{mycolor1}{rgb}{0.00000,0.44700,0.74100}%
\definecolor{mycolor2}{rgb}{0.85000,0.32500,0.09800}%
\definecolor{mycolor3}{rgb}{0.92900,0.69400,0.12500}%
\definecolor{mycolor4}{rgb}{0.49400,0.18400,0.55600}%

\pgfplotstableread[col sep=comma]{figs/offAt_mode3.csv}{\offAtIII}

\begin{tikzpicture}[trim axis left, trim axis right]
  \begin{axis}[%
    xmode=log,
    xmin=10,
    xmax=1000,
    xminorticks=true,
    ymode=log,
    ymin=5e-08,
    ymax=1e-05,
    yminorticks=true,
    xlabel={$n$},
    ylabel={$\off(\At)/\fnorm{A}$}
    ]
    \addplot [color=mycolor1, only marks,
    mark size=3.0pt, mark=triangle, mark options={solid, mycolor1}]
    table[x=n,y=qr]{\offAtIII};
    
    \addplot [color=mycolor2, only marks,
    mark size=3.0pt, mark=diamond, mark options={solid, mycolor2}]
    table[x=n,y=mgs]{\offAtIII};
    
    \addplot [color=mycolor3, only marks,
    mark size=3.0pt, mark=x, mark options={solid, mycolor3}]
    table[x=n,y=ns]{\offAtIII};
    
    \addplot [color=mycolor4, only marks,
    mark size=3.0pt, mark=star, mark options={solid, mycolor4}]
    table[x=n,y=tri]{\offAtIII};
    
    \addplot [color=black, dashed, line width=2.pt]
    table[x=n,y=bound]{\offAtIII};
  \end{axis}
\end{tikzpicture}%


%% file: figs/offAt_mode4.tex
%
\definecolor{mycolor1}{rgb}{0.00000,0.44700,0.74100}%
\definecolor{mycolor2}{rgb}{0.85000,0.32500,0.09800}%
\definecolor{mycolor3}{rgb}{0.92900,0.69400,0.12500}%
\definecolor{mycolor4}{rgb}{0.49400,0.18400,0.55600}%

\pgfplotstableread[col sep=comma]{figs/offAt_mode4.csv}{\offAtIV}

\begin{tikzpicture}[trim axis left, trim axis right]
  
  \begin{axis}[%
    xmode=log,
    xmin=10,
    xmax=1000,
    xminorticks=true,
    ymode=log,
    ymin=5e-08,
    ymax=1e-05,
    yminorticks=true,
    xlabel={$n$},
    ]
    \addplot [color=mycolor1, only marks, mark size=3.0pt,
    mark=triangle, mark options={solid, mycolor1}]
    table[x=n,y=qr]{\offAtIV};
    \addplot [color=mycolor2, only marks, mark size=3.0pt,
    mark=diamond, mark options={solid, mycolor2}]
    table[x=n,y=mgs]{\offAtIV};
    \addplot [color=mycolor3, only marks, mark size=3.0pt,
    mark=x, mark options={solid, mycolor3}]
    table[x=n,y=ns]{\offAtIV};
    \addplot [color=mycolor4, only marks, mark size=3.0pt, 
    mark=star, mark options={solid, mycolor4}]
    table[x=n,y=tri]{\offAtIV};
    \addplot [color=black, dashed, line width=2.0pt]
    table[x=n,y=bound]{\offAtIV};
  \end{axis}
  
\end{tikzpicture}%


%% file: figs/offAt_mode5.tex
%
\definecolor{mycolor1}{rgb}{0.00000,0.44700,0.74100}%
\definecolor{mycolor2}{rgb}{0.85000,0.32500,0.09800}%
\definecolor{mycolor3}{rgb}{0.92900,0.69400,0.12500}%
\definecolor{mycolor4}{rgb}{0.49400,0.18400,0.55600}%

\pgfplotstableread[col sep=comma]{figs/offAt_mode5.csv}{\offAtV}

\begin{tikzpicture}[trim axis left,trim axis right]

  \begin{axis}[
    xmode=log,
    xmin=10,
    xmax=1000,
    xminorticks=true,
    ymode=log,
    ymin=5e-08,
    ymax=1e-05,
    yminorticks=true,
    xlabel={$n$},
    ylabel={$\off(\At)/\fnorm{A}$}
    ]
    \addplot [color=mycolor1, only marks, mark size=3.0pt,
    mark=triangle, mark options={solid, mycolor1}]
    table[x=n,y=qr]{\offAtV};
    \addplot [color=mycolor2, only marks, mark size=3.0pt, 
    mark=diamond, mark options={solid, mycolor2}]
    table[x=n,y=mgs]{\offAtV};
    \addplot [color=mycolor3, only marks, mark size=3.0pt,
    mark=x, mark options={solid, mycolor3}, forget plot]
    table[x=n,y=ns]{\offAtV};
    \addplot [color=mycolor4, only marks, mark size=3.0pt,
    mark=star, mark options={solid, mycolor4}]
    table[x=n,y=tri]{\offAtV};
    \addplot [color=black, dashed, line width=2.0pt]
    table[x=n,y=bound]{\offAtV};
  \end{axis}
\end{tikzpicture}%


%% file: figs/offAt_legend.tex
\definecolor{mycolor1}{rgb}{0.00000,0.44700,0.74100}%
\definecolor{mycolor2}{rgb}{0.85000,0.32500,0.09800}%
\definecolor{mycolor3}{rgb}{0.92900,0.69400,0.12500}%
\definecolor{mycolor4}{rgb}{0.49400,0.18400,0.55600}%

\pgfplotscreateplotcyclelist{sincos}{
  color=mycolor1, only marks, mark size=3.0pt,
  mark=triangle, mark options={solid, mycolor1}\\
  color=mycolor2, only marks, mark size=3.0pt,
  mark=diamond, mark options={solid, mycolor2}\\
  color=mycolor3, only marks, mark size=3.0pt, mark=x,
  mark options={solid, mycolor3}\\
  color=mycolor4, only marks, mark size=3.0pt, mark=star,
  mark options={solid, mycolor4}\\
  color=black, dashed, line width=2.0pt\\}

\begin{tikzpicture}[trim axis left, trim axis right]
\begin{axis}[
  title = {},
  legend columns=1,
  scale only axis,
  width=0.45\textwidth,
  height=0.3\textwidth,
  xmin=0,
  xmax=1,
  ymin=0,
  ymax=1,
  hide axis,
  legend style={at={(0.575,0.6)}},
  cycle list name = sincos
  ]
  \addplot[color=mycolor1, only marks, mark size=3pt,
  mark=triangle, mark options={solid, mycolor1}](-1e10,-1e10);
  \addplot (-1e10,-1e10);
  \addplot (-1e10,-1e10);
  \addplot (-1e10,-1e10);
  \addplot (-1e10,-1e10);
  
  \legend{
  Algorithm~\ref{alg.orth-preconditioner} with HHQR,
  Algorithm~\ref{alg.orth-preconditioner} with MGS,
  Algorithm~\ref{alg.orth-preconditioner} with NS,
  Algorithm~\ref{alg.modified-tridiag}, $5n^{1/2}\ul$ };
\end{axis}
\end{tikzpicture}


%% file: figs/varykappa_mode1.tex
%
\definecolor{mycolor1}{rgb}{0.00000,0.44700,0.74100}%
\definecolor{mycolor2}{rgb}{0.85000,0.32500,0.09800}%
\definecolor{mycolor3}{rgb}{0.92900,0.69400,0.12500}%
\definecolor{mycolor4}{rgb}{0.49400,0.18400,0.55600}%

\pgfplotstableread[col sep=comma]{figs/varykappa_mode1.csv}{\varykmI}

\begin{tikzpicture}[trim axis left,trim axis right]
  \begin{axis}[%
    xmode=log,
    xmin=10,
    xmax=1e16,
    ymode=log,
    ymin=1e-16,
    ymax=150,
    ytick distance = 1e4,
    xtick distance = 1e4,
    xminorticks = true,
    yminorticks = true,
    xlabel={$\kappa_{2}(A)$},
    ylabel={$\max_{k}\ferrk$},
    ]
    \addplot[color=mycolor1, only marks,
    mark size=3.0pt, mark=triangle, mark options={solid, mycolor1}]
    table [x=kappa, y=precj] {\varykmI};
    \addplot [color=mycolor2,only marks,
    mark size=3.0pt, mark=diamond, mark options={solid, mycolor2}]
    table [x=kappa, y=uhprecj] {\varykmI};
    \addplot [color=mycolor3, only marks,
    mark size=3.0pt, mark=x, mark options={solid, mycolor3}]
    table [x=kappa, y=jacobi] {\varykmI};
    \addplot [color=mycolor4, only marks,
    mark size=3.0pt, mark=star, mark options={solid, mycolor4}]
    table [x=kappa, y=matlab] {\varykmI};
    \addplot [color=black, dashed, line width=2.0pt]
    table [x=kappa, y=bound] {\varykmI};
  \end{axis}
\end{tikzpicture}%


%% file: figs/varykappa_mode2.tex
\definecolor{mycolor1}{rgb}{0.00000,0.44700,0.74100}%
\definecolor{mycolor2}{rgb}{0.85000,0.32500,0.09800}%
\definecolor{mycolor3}{rgb}{0.92900,0.69400,0.12500}%
\definecolor{mycolor4}{rgb}{0.49400,0.18400,0.55600}%

\pgfplotstableread[col sep=comma]{figs/varykappa_mode2.csv}{\varykmII}

\begin{tikzpicture}[trim axis left,trim axis right]
  \begin{axis}[%
    xmode=log,
    xmin=10,
    xmax=1e16,
    xminorticks=true,
    ymode=log,
    ymin=1e-16,
    ymax=10,
    yminorticks=true,
    ytick distance = 1e4,
    xtick distance = 1e4,
    xlabel={$\kappa_{2}(A)$}
    ]
    \addplot[color=mycolor1, only marks,
    mark size=3.0pt, mark=triangle, mark options={solid, mycolor1}]
    table [x=kappa, y=precj] {\varykmII};
    \addplot [color=mycolor2,only marks,
    mark size=3.0pt, mark=diamond, mark options={solid, mycolor2}]
    table [x=kappa, y=uhprecj] {\varykmII};
    \addplot [color=mycolor3, only marks,
    mark size=3.0pt, mark=x, mark options={solid, mycolor3}]
    table [x=kappa, y=jacobi] {\varykmII};
    \addplot [color=mycolor4, only marks,
    mark size=3.0pt, mark=star, mark options={solid, mycolor4}]
    table [x=kappa, y=matlab] {\varykmII};
    \addplot [color=black, dashed, line width=2.0pt]
    table [x=kappa, y=bound] {\varykmII};
  \end{axis}
\end{tikzpicture}%


%% file: figs/varykappa_mode3.tex
\definecolor{mycolor1}{rgb}{0.00000,0.44700,0.74100}%
\definecolor{mycolor2}{rgb}{0.85000,0.32500,0.09800}%
\definecolor{mycolor3}{rgb}{0.92900,0.69400,0.12500}%
\definecolor{mycolor4}{rgb}{0.49400,0.18400,0.55600}%

\pgfplotstableread[col sep=comma]{figs/varykappa_mode3.csv}{\varykmIII}

\begin{tikzpicture}[trim axis left,trim axis right]
  \begin{axis}[%
    xmode=log,
    xmin=10,
    xmax=1e16,
    xminorticks=true,
    ymode=log,
    ymin=1e-15,
    ymax=1e-01,
    yminorticks=true,
    xlabel={$\kappa_{2}(A)$},
    ylabel={$\max_{k}\ferrk$},
    ytick distance = 1e4,
    xtick distance = 1e4,
    ]
    \addplot[color=mycolor1, only marks,
    mark size=3.0pt, mark=triangle, mark options={solid, mycolor1}]
    table [x=kappa, y=precj] {\varykmIII};
    \addplot [color=mycolor2,only marks,
    mark size=3.0pt, mark=diamond, mark options={solid, mycolor2}]
    table [x=kappa, y=uhprecj] {\varykmIII};
    \addplot [color=mycolor3, only marks,
    mark size=3.0pt, mark=x, mark options={solid, mycolor3}]
    table [x=kappa, y=jacobi] {\varykmIII};
    \addplot [color=mycolor4, only marks,
    mark size=3.0pt, mark=star, mark options={solid, mycolor4}]
    table [x=kappa, y=matlab] {\varykmIII};
    \addplot [color=black, dashed, line width=2.0pt]
    table [x=kappa, y=bound] {\varykmIII};
  \end{axis}
\end{tikzpicture}%


%% file: figs/varykappa_mode4.tex
\definecolor{mycolor1}{rgb}{0.00000,0.44700,0.74100}%
\definecolor{mycolor2}{rgb}{0.85000,0.32500,0.09800}%
\definecolor{mycolor3}{rgb}{0.92900,0.69400,0.12500}%
\definecolor{mycolor4}{rgb}{0.49400,0.18400,0.55600}%

\pgfplotstableread[col sep=comma]{figs/varykappa_mode4.csv}{\varykmIIII}

\begin{tikzpicture}[trim axis left,trim axis right]
  \begin{axis}[%
    xmode=log,
    xmin=10,
    xmax=1e16,
    xminorticks=true,
    ymode=log,
    ymin=1e-15,
    ymax=50,
    yminorticks=true,
    xlabel={$\kappa_{2}(A)$},
    ytick distance = 1e4,
    xtick distance = 1e4,
    ]
    \addplot[color=mycolor1, only marks,
    mark size=3.0pt, mark=triangle, mark options={solid, mycolor1}]
    table [x=kappa, y=precj] {\varykmIIII};
    \addplot [color=mycolor2,only marks,
    mark size=3.0pt, mark=diamond, mark options={solid, mycolor2}]
    table [x=kappa, y=uhprecj] {\varykmIIII};
    \addplot [color=mycolor3, only marks,
    mark size=3.0pt, mark=x, mark options={solid, mycolor3}]
    table [x=kappa, y=jacobi] {\varykmIIII};
    \addplot [color=mycolor4, only marks,
    mark size=3.0pt, mark=star, mark options={solid, mycolor4}]
    table [x=kappa, y=matlab] {\varykmIIII};
    \addplot [color=black, dashed, line width=2.0pt]
    table [x=kappa, y=bound] {\varykmIIII};
  \end{axis}
\end{tikzpicture}%


%% file: figs/varykappa_mode5.tex
\definecolor{mycolor1}{rgb}{0.00000,0.44700,0.74100}%
\definecolor{mycolor2}{rgb}{0.85000,0.32500,0.09800}%
\definecolor{mycolor3}{rgb}{0.92900,0.69400,0.12500}%
\definecolor{mycolor4}{rgb}{0.49400,0.18400,0.55600}%

\pgfplotstableread[col sep=comma]{figs/varykappa_mode5.csv}{\varykmV}

\begin{tikzpicture}[trim axis left,trim axis right]
  \begin{axis}[%
    xmode=log,
    xmin=10,
    xmax=1e16,
    xminorticks=true,
    ymode=log,
    ymin=1e-15,
    ymax=1e-01,
    yminorticks=true,
    xlabel={$\kappa_{2}(A)$},
    ylabel={$\max_{k}\ferrk$},
    ytick distance = 1e4,
    xtick distance = 1e4,
    ]
    \addplot[color=mycolor1, only marks,
    mark size=3.0pt, mark=triangle, mark options={solid, mycolor1}]
    table [x=kappa, y=precj] {\varykmV};
    \addplot [color=mycolor2,only marks,
    mark size=3.0pt, mark=diamond, mark options={solid, mycolor2}]
    table [x=kappa, y=uhprecj] {\varykmV};
    \addplot [color=mycolor3, only marks,
    mark size=3.0pt, mark=x, mark options={solid, mycolor3}]
    table [x=kappa, y=jacobi] {\varykmV};
    \addplot [color=mycolor4, only marks,
    mark size=3.0pt, mark=star, mark options={solid, mycolor4}]
    table [x=kappa, y=matlab] {\varykmV};
    \addplot [color=black, dashed, line width=2.0pt]
    table [x=kappa, y=bound] {\varykmV};
  \end{axis}
\end{tikzpicture}%


%% file: figs/varykappa_legend.tex
\definecolor{mycolor1}{rgb}{0.00000,0.44700,0.74100}%
\definecolor{mycolor2}{rgb}{0.85000,0.32500,0.09800}%
\definecolor{mycolor3}{rgb}{0.92900,0.69400,0.12500}%
\definecolor{mycolor4}{rgb}{0.49400,0.18400,0.55600}%

\pgfplotscreateplotcyclelist{sincos}{
  color=mycolor1, only marks, mark size=3.0pt,
  mark=triangle, mark options={solid, mycolor1}\\
  color=mycolor2,  only marks, mark size=3.0pt,
  mark=diamond, mark options={solid, mycolor2}\\
  color=mycolor3, only marks, mark size=3.0pt, mark=x,
  mark options={solid, mycolor3}\\
  color=mycolor4, only marks, mark size=3.0pt, mark=star,
  mark options={solid, mycolor4}\\
  color=black, dashed, line width=2.0pt\\}

\begin{tikzpicture}[trim axis left, trim axis right]
\begin{axis}[
  title = {},
  legend columns=1,
  scale only axis,
  width=0.45\textwidth,
  height=0.3\textwidth,
  xmin=0,
  xmax=1,
  ymin=0,
  ymax=1,
  hide axis,
  legend style={at={(0.6,0.6)}},
  cycle list name = sincos
  ]
  \addplot (-1e10,-1e10);
  \addplot (-1e10,-1e10);
  \addplot (-1e10,-1e10);
  \addplot (-1e10,-1e10);
  \addplot (-1e10,-1e10);
  
  \legend{{MP2Jacobi}, {MP3Jacobi}, {Jacobi},
           {MATLAB}, $7n\scond(\At)u$};
\end{axis}
\end{tikzpicture}


%% file: figs/varydim_mode1.tex
%
\definecolor{mycolor1}{rgb}{0.00000,0.44700,0.74100}%
\definecolor{mycolor2}{rgb}{0.85000,0.32500,0.09800}%
\definecolor{mycolor3}{rgb}{0.92900,0.69400,0.12500}%
\definecolor{mycolor4}{rgb}{0.49400,0.18400,0.55600}%

\pgfplotstableread[col sep=comma]{figs/varydim_mode1.csv}{\varydimI}

\begin{tikzpicture}[trim axis left,trim axis right]
  \begin{axis}[%
    xmode=log,
    xmin=10,
    xmax=1e3,
    xminorticks=true,
    ymode=log,
    ymin=1e-16,
    ymax=1e-4,
    yminorticks=true,
    xlabel={$n$},
    ylabel={$\max_{k}\ferrk$},
    ytick distance = 1e4,
    ]
    \addplot[color=mycolor1, only marks,
    mark size=3.0pt, mark=triangle, mark options={solid, mycolor1}]
    table [x=kappa, y=precj] {\varydimI};
    \addplot [color=mycolor2,only marks,
    mark size=3.0pt, mark=diamond, mark options={solid, mycolor2}]
    table [x=kappa, y=uhprecj] {\varydimI};
    \addplot [color=mycolor3, only marks,
    mark size=3.0pt, mark=x, mark options={solid, mycolor3}]
    table [x=kappa, y=jacobi] {\varydimI};
    \addplot [color=mycolor4, only marks,
    mark size=3.0pt, mark=star, mark options={solid, mycolor4}]
    table [x=kappa, y=matlab] {\varydimI};
    \addplot [color=black, dashed, line width=2.0pt]
    table [x=kappa, y=bound] {\varydimI};
  \end{axis}
\end{tikzpicture}%


%% file: figs/varydim_mode2.tex
%
\definecolor{mycolor1}{rgb}{0.00000,0.44700,0.74100}%
\definecolor{mycolor2}{rgb}{0.85000,0.32500,0.09800}%
\definecolor{mycolor3}{rgb}{0.92900,0.69400,0.12500}%
\definecolor{mycolor4}{rgb}{0.49400,0.18400,0.55600}%

\pgfplotstableread[col sep=comma]{figs/varydim_mode2.csv}{\varydimII}

\begin{tikzpicture}[trim axis left,trim axis right]
  \begin{axis}[%
    xmode=log,
    xmin=10,
    xmax=1e3,
    xminorticks=true,
    ymode=log,
    ymin=1e-16,
    ymax=1e-5,
    yminorticks=true,
    xlabel={$n$},
    ytick distance = 1e4,
    ]
    \addplot[color=mycolor1, only marks,
    mark size=3.0pt, mark=triangle, mark options={solid, mycolor1}]
    table [x=kappa, y=precj] {\varydimII};
    \addplot [color=mycolor2,only marks,
    mark size=3.0pt, mark=diamond, mark options={solid, mycolor2}]
    table [x=kappa, y=uhprecj] {\varydimII};
    \addplot [color=mycolor3, only marks,
    mark size=3.0pt, mark=x, mark options={solid, mycolor3}]
    table [x=kappa, y=jacobi] {\varydimII};
    \addplot [color=mycolor4, only marks,
    mark size=3.0pt, mark=star, mark options={solid, mycolor4}]
    table [x=kappa, y=matlab] {\varydimII};
    \addplot [color=black, dashed, line width=2.0pt]
    table [x=kappa, y=bound] {\varydimII};
  \end{axis}
\end{tikzpicture}%


%% file: figs/varydim_mode3.tex
%
\definecolor{mycolor1}{rgb}{0.00000,0.44700,0.74100}%
\definecolor{mycolor2}{rgb}{0.85000,0.32500,0.09800}%
\definecolor{mycolor3}{rgb}{0.92900,0.69400,0.12500}%
\definecolor{mycolor4}{rgb}{0.49400,0.18400,0.55600}%

\pgfplotstableread[col sep=comma]{figs/varydim_mode3.csv}{\varydimIII}

\begin{tikzpicture}[trim axis left,trim axis right]
  \begin{axis}[%
    xmode=log,
    xmin=10,
    xmax=1e3,
    xminorticks=true,
    ymode=log,
    ymin=1e-16,
    ymax=1e-5,
    yminorticks=true,
    ylabel={$\max_{k}\ferrk$},
    xlabel={$n$},
    ytick distance = 1e4,
    ]
    \addplot[color=mycolor1, only marks,
    mark size=3.0pt, mark=triangle, mark options={solid, mycolor1}]
    table [x=kappa, y=precj] {\varydimIII};
    \addplot [color=mycolor2,only marks,
    mark size=3.0pt, mark=diamond, mark options={solid, mycolor2}]
    table [x=kappa, y=uhprecj] {\varydimIII};
    \addplot [color=mycolor3, only marks,
    mark size=3.0pt, mark=x, mark options={solid, mycolor3}]
    table [x=kappa, y=jacobi] {\varydimIII};
    \addplot [color=mycolor4, only marks,
    mark size=3.0pt, mark=star, mark options={solid, mycolor4}]
    table [x=kappa, y=matlab] {\varydimIII};
    \addplot [color=black, dashed, line width=2.0pt]
    table [x=kappa, y=bound] {\varydimIII};
  \end{axis}
\end{tikzpicture}%


%% file: figs/varydim_mode4.tex
%
\definecolor{mycolor1}{rgb}{0.00000,0.44700,0.74100}%
\definecolor{mycolor2}{rgb}{0.85000,0.32500,0.09800}%
\definecolor{mycolor3}{rgb}{0.92900,0.69400,0.12500}%
\definecolor{mycolor4}{rgb}{0.49400,0.18400,0.55600}%

\pgfplotstableread[col sep=comma]{figs/varydim_mode4.csv}{\varydimIV}

\begin{tikzpicture}[trim axis left,trim axis right]
  \begin{axis}[%
    xmode=log,
    xmin=10,
    xmax=1e3,
    xminorticks=true,
    ymode=log,
    ymin=1e-16,
    ymax=1e-4,
    yminorticks=true,
    xlabel={$n$},
    ytick distance = 1e4,
    ]
    \addplot[color=mycolor1, only marks,
    mark size=3.0pt, mark=triangle, mark options={solid, mycolor1}]
    table [x=kappa, y=precj] {\varydimIV};
    \addplot [color=mycolor2,only marks,
    mark size=3.0pt, mark=diamond, mark options={solid, mycolor2}]
    table [x=kappa, y=uhprecj] {\varydimIV};
    \addplot [color=mycolor3, only marks,
    mark size=3.0pt, mark=x, mark options={solid, mycolor3}]
    table [x=kappa, y=jacobi] {\varydimIV};
    \addplot [color=mycolor4, only marks,
    mark size=3.0pt, mark=star, mark options={solid, mycolor4}]
    table [x=kappa, y=matlab] {\varydimIV};
    \addplot [color=black, dashed, line width=2.0pt]
    table [x=kappa, y=bound] {\varydimIV};
  \end{axis}
\end{tikzpicture}%


%% file: figs/varydim_mode5.tex
%
\definecolor{mycolor1}{rgb}{0.00000,0.44700,0.74100}%
\definecolor{mycolor2}{rgb}{0.85000,0.32500,0.09800}%
\definecolor{mycolor3}{rgb}{0.92900,0.69400,0.12500}%
\definecolor{mycolor4}{rgb}{0.49400,0.18400,0.55600}%

\pgfplotstableread[col sep=comma]{figs/varydim_mode5.csv}{\varydimV}

\begin{tikzpicture}[trim axis left,trim axis right]
  \begin{axis}[%
    xmode=log,
    xmin=10,
    xmax=1e3,
    xminorticks=true,
    ymode=log,
    ymin=1e-16,
    ymax=1e-6,
    yminorticks=true,
    xlabel={$n$},
    ylabel={$\max_{k}\ferrk$},
    ytick distance = 1e4,
    ]
    \addplot[color=mycolor1, only marks,
    mark size=3.0pt, mark=triangle, mark options={solid, mycolor1}]
    table [x=kappa, y=precj] {\varydimV};
    \addplot [color=mycolor2,only marks,
    mark size=3.0pt, mark=diamond, mark options={solid, mycolor2}]
    table [x=kappa, y=uhprecj] {\varydimV};
    \addplot [color=mycolor3, only marks,
    mark size=3.0pt, mark=x, mark options={solid, mycolor3}]
    table [x=kappa, y=jacobi] {\varydimV};
    \addplot [color=mycolor4, only marks,
    mark size=3.0pt, mark=star, mark options={solid, mycolor4}]
    table [x=kappa, y=matlab] {\varydimV};
    \addplot [color=black, dashed, line width=2.0pt]
    table [x=kappa, y=bound] {\varydimV};
  \end{axis}
\end{tikzpicture}%


%% file: figs/varydim_legend.tex
\definecolor{mycolor1}{rgb}{0.00000,0.44700,0.74100}%
\definecolor{mycolor2}{rgb}{0.85000,0.32500,0.09800}%
\definecolor{mycolor3}{rgb}{0.92900,0.69400,0.12500}%
\definecolor{mycolor4}{rgb}{0.49400,0.18400,0.55600}%

\pgfplotscreateplotcyclelist{sincos}{
  color=mycolor1, only marks, mark size=3.0pt,
  mark=triangle, mark options={solid, mycolor1}\\
  color=mycolor2, only marks, mark size=3.0pt,
  mark=diamond, mark options={solid, mycolor2}\\
  color=mycolor3, only marks, mark size=3.0pt, 
  mark=x, mark options={solid, mycolor3}\\
  color=mycolor4, only marks, mark size=3.0pt, mark=star,
  mark options={solid, mycolor4}\\
  color=black, dashed, line width=2.0pt\\}

\begin{tikzpicture}[trim axis left, trim axis right]
\begin{axis}[
  title = {},
  legend columns=1,
  scale only axis,
  width=0.45\textwidth,
  height=0.3\textwidth,
  xmin=0,
  xmax=1,
  ymin=0,
  ymax=1,
  hide axis,
  legend style={at={(0.6,0.6)}},
  cycle list name = sincos
  ]
  \addplot (-1e10,-1e10);
  \addplot (-1e10,-1e10);
  \addplot (-1e10,-1e10);
  \addplot (-1e10,-1e10);
  \addplot (-1e10,-1e10);
  
  \legend{{MP2Jacobi}, {MP3Jacobi}, {Jacobi},
           {MATLAB}, $7n\scond(\At)u$};
\end{axis}
\end{tikzpicture}


%% file: figs/rand_matrix.tex
%
\definecolor{mycolor1}{rgb}{0.00000,0.44700,0.74100}%
\definecolor{mycolor2}{rgb}{0.85000,0.32500,0.09800}%
\definecolor{mycolor3}{rgb}{0.92900,0.69400,0.12500}%
\definecolor{mycolor4}{rgb}{0.49400,0.18400,0.55600}%

\pgfplotstableread[col sep=comma]{figs/random_matrix.csv}{\random}

\begin{tikzpicture}[trim axis left, trim axis right]
  \begin{axis}[%
    width = 0.875\textwidth,
    height = 0.325\textwidth,
    ymode=log,
    xmin=1,
    xmax=100,
    xminorticks=true,
    ymode=log,
    ymin=1e-17,
    ymax=1e-8,
    yminorticks=true,
    xlabel={$k$th largest eigenvalue},
    ylabel={$\ferrk$},
    ytick distance = 1e3,
    ]
    \addplot[color=mycolor1, only marks,
    mark size=3.0pt, mark=triangle, mark options={solid, mycolor1}]
    table [x=k, y=mp2] {\random};
    \addplot [color=mycolor2, only marks,
    mark size=3.0pt, mark=diamond, mark options={solid, mycolor2}]
    table [x=k, y=mp3] {\random};
    \addplot [color=mycolor3, only marks,
    mark size=3.0pt, mark=x, mark options={solid, mycolor3}]
    table [x=k, y=jcb] {\random};
    \addplot [color=mycolor4, only marks,
    mark size=3.0pt, mark=star, mark options={solid, mycolor4}]
    table [x=k, y=matlab] {\random};
    \addplot [color=black, dashed, line width=2.0pt]
    table [x=k, y=bound] {\random};
  \end{axis}
\end{tikzpicture}%


%% file: figs/whisky.tex
%
\definecolor{mycolor1}{rgb}{0.00000,0.44700,0.74100}%
\definecolor{mycolor2}{rgb}{0.85000,0.32500,0.09800}%
\definecolor{mycolor3}{rgb}{0.92900,0.69400,0.12500}%
\definecolor{mycolor4}{rgb}{0.49400,0.18400,0.55600}%

\pgfplotstableread[col sep=comma]{figs/whiskycorr.csv}{\whisky}

\begin{tikzpicture}[trim axis left, trim axis right]
  \begin{axis}[%
    width = 0.875\textwidth,
    height = 0.325\textwidth,
    ymode=log,
    xmin=1,
    xmax=86,
    xminorticks=true,
    ymode=log,
    ymin=1e-18,
    ymax=1e2,
    yminorticks=true,
    xlabel={$k$th largest eigenvalue},
    ylabel={$\ferrk$},
    ytick distance = 1e4,
    ]
    \addplot[color=mycolor1, only marks,
    mark size=3.0pt, mark=triangle, mark options={solid, mycolor1}]
    table [x=k, y=mp2] {\whisky};
    \addplot [color=mycolor2, only marks,
    mark size=3.0pt, mark=diamond, mark options={solid, mycolor2}]
    table [x=k, y=mp3] {\whisky};
    \addplot [color=mycolor3, only marks,
    mark size=3.0pt, mark=x, mark options={solid, mycolor3}]
    table [x=k, y=jcb] {\whisky};
    \addplot [color=mycolor4, only marks,
    mark size=3.0pt, mark=star, mark options={solid, mycolor4}]
    table [x=k, y=matlab] {\whisky};
    \addplot [color=black, dashed, line width=2.0pt]
    table [x=k, y=bound] {\whisky};
  \end{axis}
\end{tikzpicture}%


%% file: figs/covariance.tex
\definecolor{mycolor1}{rgb}{0.00000,0.44700,0.74100}%
\definecolor{mycolor2}{rgb}{0.85000,0.32500,0.09800}%
\definecolor{mycolor3}{rgb}{0.92900,0.69400,0.12500}%
\definecolor{mycolor4}{rgb}{0.49400,0.18400,0.55600}%

\pgfplotstableread[col sep=comma]{figs/covariance.csv}{\cov}

\begin{tikzpicture}[trim axis left, trim axis right]
  \begin{axis}[%
    width = 0.875\textwidth,
    height = 0.325\textwidth,
    ymode=log,
    xmin=1,
    xmax=82,
    xminorticks=true,
    ymode=log,
    ymin=1e-17,
    ymax=1e-5,
    yminorticks=true,
    xlabel={$k$th largest eigenvalue},
    ylabel={$\ferrk$},
    ytick distance = 1e2,
    xtick distance = 10,
    ]
    \addplot[color=mycolor1, only marks,
    mark size=3.0pt, mark=triangle, mark options={solid, mycolor1}]
    table [x=k, y=mp2] {\cov};
    \addplot [color=mycolor2, only marks,
    mark size=3.0pt, mark=diamond, mark options={solid, mycolor2}]
    table [x=k, y=mp3] {\cov};
    \addplot [color=mycolor3, only marks,
    mark size=3.0pt, mark=x, mark options={solid, mycolor3}]
    table [x=k, y=jcb] {\cov};
    \addplot [color=mycolor4, only marks,
    mark size=3.0pt, mark=star, mark options={solid, mycolor4}]
    table [x=k, y=matlab] {\cov};
    \addplot [color=black, dashed, line width=2.0pt]
    table [x=k, y=bound] {\cov};
  \end{axis}
\end{tikzpicture}%

%% file: figs/nos7.tex
\definecolor{mycolor1}{rgb}{0.00000,0.44700,0.74100}%
\definecolor{mycolor2}{rgb}{0.85000,0.32500,0.09800}%
\definecolor{mycolor3}{rgb}{0.92900,0.69400,0.12500}%
\definecolor{mycolor4}{rgb}{0.49400,0.18400,0.55600}%

\pgfplotstableread[col sep=comma]{figs/lau.csv}{\lau}

\begin{tikzpicture}[trim axis left, trim axis right]
  \begin{axis}[%
    width = 0.875\textwidth,
    height = 0.325\textwidth,
    ymode=log,
    xmin=1,
    xmax=500,
    xminorticks=true,
    ymode=log,
    ymin=1e-17,
    ymax=1e-5,
    yminorticks=true,
    xlabel={$k$th largest eigenvalue},
    ylabel={$\ferrk$},
    ytick distance = 1e2,
    ]
    \addplot[color=mycolor1, only marks,
    mark size=3.0pt, mark=triangle, mark options={solid, mycolor1}]
    table [x=k, y=mp2] {\lau};
    \addplot [color=mycolor2, only marks,
    mark size=3.0pt, mark=diamond, mark options={solid, mycolor2}]
    table [x=k, y=mp3] {\lau};
    \addplot [color=mycolor3, only marks,
    mark size=3.0pt, mark=x, mark options={solid, mycolor3}]
    table [x=k, y=jcb] {\lau};
    \addplot [color=mycolor4, only marks,
    mark size=3.0pt, mark=star, mark options={solid, mycolor4}]
    table [x=k, y=matlab] {\lau};
    \addplot [color=black, dashed, line width=2.0pt]
    table [x=k, y=bound] {\lau};
  \end{axis}
\end{tikzpicture}%

%% file: figs/special_legend.tex
\definecolor{mycolor1}{rgb}{0.00000,0.44700,0.74100}%
\definecolor{mycolor2}{rgb}{0.85000,0.32500,0.09800}%
\definecolor{mycolor3}{rgb}{0.92900,0.69400,0.12500}%
\definecolor{mycolor4}{rgb}{0.49400,0.18400,0.55600}%

\pgfplotscreateplotcyclelist{sincos}{
  color=mycolor1, only marks, mark size=3.0pt,
  mark=triangle, mark options={solid, mycolor1}\\
  color=mycolor2, only marks, mark size=3.0pt,
  mark=diamond, mark options={solid, mycolor2}\\
  color=mycolor3, only marks, mark size=3.0pt, mark=x,
  mark options={solid, mycolor3}\\
  color=mycolor4, only marks, mark size=3.0pt, mark=star,
  mark options={solid, mycolor4}\\
  color=black, dashed, line width=2.0pt\\}

\vspace{.3cm}
\begin{tikzpicture}[trim axis left, trim axis right]
\begin{axis}[
  title = {},
  legend columns=5,
  scale only axis,
  width=0.4\textwidth,
  xmin=0,
  xmax=1,
  ymin=0,
  ymax=1,
  hide axis,
  /tikz/every even column/.append style={column sep=0.6cm},
  legend style={at={(0.5,1)}},
  cycle list name = sincos
  ]
  \addplot (-1e10,-1e10);
  \addplot (-1e10,-1e10);
  \addplot (-1e10,-1e10);
  \addplot (-1e10,-1e10);
  \addplot (-1e10,-1e10);
  
  \legend{{MP2Jacobi},
  {MP3Jacobi},
  {Jacobi},
  {MATLAB},
  $7n\scond(\At)u$};
\end{axis}
\end{tikzpicture}

\vspace{-4.1cm}


%% file: figs/reducekappa_mode1.tex
%
\definecolor{mycolor1}{rgb}{0.00000,0.44700,0.74100}%
\definecolor{mycolor2}{rgb}{0.85000,0.32500,0.09800}%
\definecolor{mycolor3}{rgb}{0.92900,0.69400,0.12500}%
\definecolor{mycolor4}{rgb}{0.49400,0.18400,0.55600}%

\pgfplotstableread[col sep=comma]{figs/reducekappa_mode1.csv}{\rdkI}

\begin{tikzpicture}[trim axis left,trim axis right]
  \begin{axis}[%
    xmode=log,
    xmin=9.5,
    xmax=1e12,
    xminorticks=true,
    ymode=log,
    ymin=1e-12,
    ymax=50,
    yminorticks=true,
    xlabel={$\kappa_{2}(A)$},
    ylabel={Ratio},
    xtick distance = 1e4,
    ytick distance = 1e4,
    ]
    \addplot[color=mycolor1, only marks,
    mark size=3.0pt, mark=triangle, mark options={solid, mycolor1}]
    table [x=ka, y=DAD] {\rdkI};
    \addplot [color=mycolor2,only marks,
    mark size=3.0pt, mark=diamond, mark options={solid, mycolor2}]
    table [x=ka, y=DAtD] {\rdkI};
  \end{axis}
  
\end{tikzpicture}%


%% file: figs/reducekappa_mode2.tex
%
\definecolor{mycolor1}{rgb}{0.00000,0.44700,0.74100}%
\definecolor{mycolor2}{rgb}{0.85000,0.32500,0.09800}%
\definecolor{mycolor3}{rgb}{0.92900,0.69400,0.12500}%
\definecolor{mycolor4}{rgb}{0.49400,0.18400,0.55600}%

\pgfplotstableread[col sep=comma]{figs/reducekappa_mode2.csv}{\rdkII}

\begin{tikzpicture}[trim axis left,trim axis right]
  \begin{axis}[%
    xmode=log,
    xmin=9.5,
    xmax=1.1e12,
    xminorticks=true,
    ymode=log,
    ymin=1e-12,
    ymax=50,
    yminorticks=true,
    xlabel={$\kappa_{2}(A)$},
    xtick distance = 1e4,
    ytick distance = 1e4,
    ]
    \addplot[color=mycolor1, only marks,
    mark size=3.0pt, mark=triangle, mark options={solid, mycolor1}]
    table [x=ka, y=DAD] {\rdkII};
    \addplot [color=mycolor2,only marks,
    mark size=3.0pt, mark=diamond, mark options={solid, mycolor2}]
    table [x=ka, y=DAtD] {\rdkII};
  \end{axis}
\end{tikzpicture}%


%% file: figs/reducekappa_mode3.tex
%
\definecolor{mycolor1}{rgb}{0.00000,0.44700,0.74100}%
\definecolor{mycolor2}{rgb}{0.85000,0.32500,0.09800}%
\definecolor{mycolor3}{rgb}{0.92900,0.69400,0.12500}%
\definecolor{mycolor4}{rgb}{0.49400,0.18400,0.55600}%

\pgfplotstableread[col sep=comma]{figs/reducekappa_mode3.csv}{\rdkIII}

\begin{tikzpicture}[trim axis left,trim axis right]
  \begin{axis}[%
    xmode=log,
    xmin=9.5,
    xmax=1.1e12,
    xminorticks=true,
    ymode=log,
    ymin=1e-12,
    ymax=50,
    yminorticks=true,
    xlabel={$\kappa_{2}(A)$},
    ylabel={Ratio},
    xtick distance = 1e4,
    ytick distance = 1e4,
    ]
    \addplot[color=mycolor1, only marks,
    mark size=3.0pt, mark=triangle, mark options={solid, mycolor1}]
    table [x=ka, y=DAD] {\rdkIII};
    \addplot [color=mycolor2,only marks,
    mark size=3.0pt, mark=diamond, mark options={solid, mycolor2}]
    table [x=ka, y=DAtD] {\rdkIII};
  \end{axis}
  
\end{tikzpicture}%


%% file: figs/reducekappa_mode4.tex
%
\definecolor{mycolor1}{rgb}{0.00000,0.44700,0.74100}%
\definecolor{mycolor2}{rgb}{0.85000,0.32500,0.09800}%
\definecolor{mycolor3}{rgb}{0.92900,0.69400,0.12500}%
\definecolor{mycolor4}{rgb}{0.49400,0.18400,0.55600}%

\pgfplotstableread[col sep=comma]{figs/reducekappa_mode2.csv}{\rdkIV}

\begin{tikzpicture}[trim axis left,trim axis right]
  \begin{axis}[%
    xmode=log,
    xmin=9.5,
    xmax=1.1e12,
    xminorticks=true,
    ymode=log,
    ymin=1e-12,
    ymax=50,
    yminorticks=true,
    xlabel={$\kappa_{2}(A)$},
    xtick distance = 1e4,
    ytick distance = 1e4,
    ]
    \addplot[color=mycolor1, only marks,
    mark size=3.0pt, mark=triangle, mark options={solid, mycolor1}]
    table [x=ka, y=DAD] {\rdkIV};
    \addplot [color=mycolor2,only marks,
    mark size=3.0pt, mark=diamond, mark options={solid, mycolor2}]
    table [x=ka, y=DAtD] {\rdkIV};
  \end{axis}
\end{tikzpicture}%


%% file: figs/reducekappa_mode5.tex
%
\definecolor{mycolor1}{rgb}{0.00000,0.44700,0.74100}%
\definecolor{mycolor2}{rgb}{0.85000,0.32500,0.09800}%
\definecolor{mycolor3}{rgb}{0.92900,0.69400,0.12500}%
\definecolor{mycolor4}{rgb}{0.49400,0.18400,0.55600}%

\pgfplotstableread[col sep=comma]{figs/reducekappa_mode5.csv}{\rdkV}

\begin{tikzpicture}[trim axis left,trim axis right]
  \begin{axis}[%
    xmode=log,
    xmin=9.5,
    xmax=1.1e12,
    xminorticks=true,
    ymode=log,
    ymin=1e-12,
    ymax=50,
    yminorticks=true,
    xlabel={$\kappa_{2}(A)$},
    ylabel={Ratio},
    xtick distance = 1e4,
    ytick distance = 1e4,
    ]
    \addplot[color=mycolor1, only marks,
    mark size=3.0pt, mark=triangle, mark options={solid, mycolor1}]
    table [x=ka, y=DAD] {\rdkV};
    \addplot [color=mycolor2,only marks,
    mark size=3.0pt, mark=diamond, mark options={solid, mycolor2}]
    table [x=ka, y=DAtD] {\rdkV};
  \end{axis}
  
\end{tikzpicture}%


%% file: figs/reducekappa_legend.tex
\definecolor{mycolor1}{rgb}{0.00000,0.44700,0.74100}%
\definecolor{mycolor2}{rgb}{0.85000,0.32500,0.09800}%
\definecolor{mycolor3}{rgb}{0.92900,0.69400,0.12500}%
\definecolor{mycolor4}{rgb}{0.49400,0.18400,0.55600}%

\pgfplotscreateplotcyclelist{sincos}{
  color=mycolor1, only marks, mark size=3pt,
  mark=triangle, mark options={solid, mycolor1}\\
  color=mycolor2, only marks, mark size=3pt,
  mark=diamond, mark options={solid, mycolor2}\\}

\begin{tikzpicture}[trim axis left, trim axis right]
\begin{axis}[
  title = {},
  legend columns=1,
  scale only axis,
  width=0.45\textwidth,
  xmin=0,
  xmax=1,
  ymin=0,
  ymax=1,
  hide axis,
  legend style={at={(0.6,0.5)}},
  cycle list name = sincos
  ]
  \addplot (-1e10,-1e10);
  \addplot (-1e10,-1e10);
  \addplot (-1e10,-1e10);
  \addplot (-1e10,-1e10);
  \addplot (-1e10,-1e10);
  
  \legend{$\scond(A)/\k(A)$, $\scond(\At)/\k(A)$};
\end{axis}
\end{tikzpicture}
